\documentclass{article}
\usepackage{amssymb,amsmath,bm,geometry,graphics,graphicx,url,color}

\def\no{\noindent}
\def\pmatrix{\left(\begin{array}}
\def\endpmatrix{\end{array}\right)}

\def\RR{\mathbb{R}}

\def\H{{\cal H}}
\def\I{{\cal I}}

\def\P{{\cal P}}

\def\dd{\mathrm{d}}
\def\ee{\mathrm{e}}
\def\diag{\mathrm{diag}}

\def\ii{\mathrm{i}}

\newtheorem{theo}{Theorem}
\newtheorem{lem}{Lemma}

\newtheorem{rem}{Remark}
\newtheorem{defi}{Definition}
\newtheorem{crite}{Criterion}
\def\proof{\noindent\underline{Proof}\quad}
\def\QED{\mbox{~$\Box{~}$}}

\def\bfb{{\bm{b}}}
\def\bfc{{\bm{c}}}

\def\bfe{{\bm{e}}}

\def\bfp{{\bm{p}}}
\def\bfq{{\bm{q}}}

\def\bfu{{\bm{u}}}

\def\bfw{{\bm{w}}}

\def\bfy{{\bm{y}}}

\def\bfgamma{{\bm{\gamma}}}

\def\bfeta{{\bm{\eta}}}

\def\bfpsi{{\bm{\psi}}}

\def\blue{}

\title{On the effectiveness of spectral methods for the numerical solution of multi-frequency highly-oscillatory Hamiltonian problems}

\author{Luigi Brugnano\,\footnote{Dipartimento di Matematica e Informatica ``U.\,Dini'', Universit\`a di Firenze, Viale Morgagni 67/A, 50134 Firenze, Italy. E-mail:~{\tt luigi.brugnano@unifi.it}} 
\and Juan I.\,Montijano\thanks{Departamento  de Matem\'{a}tica Aplicada, Universidad de Zaragoza, Pza. San Francisco s/n,
50009 Zaragoza, Spain. E-mail:~{\tt monti@unizar.es}}
\and Luis R\'andez\thanks{Departamento  de Matem\'{a}tica Aplicada, Universidad de Zaragoza, Pza. San Francisco s/n,
50009 Zaragoza, Spain. E-mail:~{\tt randez@unizar.es}} 
}

\begin{document}
\maketitle

\begin{abstract} Multi-frequency, highly-oscillatory Hamiltonian problems derive from the mathematical modelling of many real life applications. We here propose a variant of Hamiltonian Boundary Value Methods (HBVMs), which is able to efficiently deal with the numerical solution of such problems.

\medskip
\no{\bf Keywords:} Multi-frequency highly-oscillatory problems, Hamiltonian problems, energy-conser\-ving methods, spectral methods, Legendre polynomials, Hamiltonian Boundary Value Methods, HBVMs.

\medskip
\no{\bf MSC:} 65P10, 65L05.

\end{abstract}

\section{Introduction}\label{intro}
Multi-frequency highly-oscillatory Hamiltonian problems appear often in mathematical models of real life applications such as molecular dynamics \cite{TuBe91}  or multibody mechanical systems \cite{ShHaJa04, SoMa2011}. They also occur when solving Hamiltonian PDEs by means of a proper space-semidiscretization. One common feature of this class of problems is that the spectrum of the Jacobian matrix of the vector field has one or more eigenvalues located on the imaginary axis and with a very large modulus.
 We here consider the efficient numerical solution of such problems. For purposes of analysis and to sketch the main facts about the methods, we consider the model problem
\begin{equation}\label{sec}
\ddot q + A^2 q+ \nabla f(q) = 0, \quad t\ge0, \qquad q(0) = q_0, ~\, \dot q(0) =\dot q_0\in \RR^m,
\end{equation}
where, without loss of generality, we  can assume $A$ to be a symmetric and positive definite (spd) matrix: in fact, possible zero eigenvalues could be, e.g., set to 1, then moving the residual to the $\nabla f(q)$ term. The equation (\ref{sec}) is defined by the separable Hamiltonian 
\begin{equation}\label{Hsec}
H(q,\dot q) = \frac{1}2(\|{\blue \dot q}\|_2^2+\|A {\blue q}\|_2^2) + f(q),
\end{equation}
with $f$ a regular enough function.  Moreover, in (\ref{sec}) we assume the nonlinear term to be ``small'' when compared with the linear part. Consequently, we shall assume
\begin{equation}\label{a2large}
\|A\| =: \omega \gg \| \nabla f\|,
\end{equation}
where the last inequality holds in a suitable domain containing the trajectory solution. Clearly, since $A$ is sdp, then there exists an orthogonal matrix $Q$ such that
\begin{equation}\label{Asdp} A = Q\Lambda Q^\top, \qquad \Lambda = \diag(\lambda_1,\dots,\lambda_m),\end{equation}
where $0<\lambda_1\le\dots\le\lambda_m\le \omega$ are the eigenvalues of $A$.
Problem (\ref{sec}) can be cast in first order form,
by setting
\begin{equation}\label{sec1}
y = \pmatrix{c} q\\ p\endpmatrix, \qquad p = A^{-1}\dot q, \qquad J_2=\pmatrix{cc} &1\\-1\endpmatrix, \qquad {\blue \tilde f(y) = \pmatrix{c} A^{-1}\nabla f(q) \\ 0\endpmatrix},
\end{equation}
and denoting ~$p_0 = A^{-1}\dot q_0$, as
\begin{equation}\label{sec2}
\dot y = J_2\otimes A\, y + J_2\otimes I_m {\blue \tilde f(y)}, \qquad y(0) = y_0 \equiv \pmatrix{c} q_0\\ p_0\endpmatrix.
\end{equation}
Consequently, the used arguments naturally extends to first order {\blue Hamiltonian} problems in the form
\begin{equation}\label{first}
\dot y = J\left[ A y + \nabla f(y)\right], \qquad y(0)=y_0\in\RR^{2m}, \qquad J = J_2\otimes I_m,
\end{equation}
with Hamiltonian
\begin{equation}\label{Hfirst}
H(y) = \frac{1}2 y^\top Ay+ f(y),
\end{equation}
and $A\in\RR^{2m\times 2m}$ a spd matrix formally still satisfying (\ref{a2large}).

The numerical solution of this class of problems presents two important difficulties. On one side, the fact that the matrix $JA$ of the linear term in (\ref{first}) has some large pure imaginary eigenvalues makes the system stiff oscillatory and  the stepsize must be small enough to guarantee that $\omega h$, with $\omega$ the largest eigenvalue modulus, belongs to the stability domain of the method.  This can be very restrictive unless the numerical method has adequate stability properties.

On the other hand, the solution $y(t)$ of (\ref{first}) can also be highly oscillatory which means that their derivatives can behave as the powers of $\omega$, that is, $y^{(j)}(t) =\mathcal{O}(\omega^j)$. Since the numerical methods are usually based on Taylor expansions and for a method of order $p$ the leading term of its local truncation error is of the order of $h^{p+1} y^{(p+1)}(t_n)$, then the error will behave as $(\omega h)^{p+1}$. Consequently, to have a small error in the numerical solution, the stepsize must again satisfy $\omega h<1$.  Otherwise, even though the error can be bounded, the error will not decrease with the stepsize $h$ until $\omega h<1$. That is, its observed numerical order can be zero  for larger values of the stepsize (see, e.g., \cite{BCMR16}).

The numerical solution of highly oscillatory problems has been the subject of many researches in the last years.
Trigonometric methods \cite{Cano2013, Cohen2006, De1979, GaSaSk98, Ga1961, HoOs10} are a class of explicit exponential methods intended for second order problems. They provide a bounded numerical solution and can be symplectic. However, even though they integrate exactly linear problems, they can present instabilities for $\omega h\ge 2\pi$ and their numerical order is zero unless $\omega h\le 1$. Functionally fitted methods  provide a generalization of the previous approach which, however, still suffers from stepsize restrictions (see, e.g., \cite[Thm.\,3.2]{LiWu2016}).

Exponential methods \cite{HoOs10} have been proved to be efficient for the solution of highly oscillatory problems coming from the semidiscretization of semilinear Hamiltonian PDEs \cite{HoOs2005, CaGo15}.  These methods require the computation of matrix exponentials, that can be expensive if the order of the method is high or if the method advances with a variable stepsize strategy. Again, to ensure the right numerical order, the stepsize must be small.

When the problem has a single high-frequency, the numerical methods can exploit such a feature for efficiently solving it.  This is the case for example
of the so called multi-revolution methods \cite{CaJaMoRa04, PeJaYe97} and the averaging or stroboscopic methods \cite{CaChMuSS11, CaChFa09, Cohen2006}.  These classes of methods combine outer integrators, that adapt to the scale of the slow components, with inner integrators, that adapt to the fast components.  They have proved to be efficient with problems that have one high frequency, but we are not aware of any result with problems with several high frequencies.
Related to this methods are multiscale techniques \cite{ArEnKiLeTs13} and parareal methods. 

A different approach is used in \cite{BCMR16} to solve second order problems with  one high frequency. These methods use a combination of Taylor and Fourier expansions to follow the high oscillations and integrate exactly linear problems, which make them stable and they can integrate with large stepsizes. Nevertheless, when multiple frequencies are present and/or they are not a priori known, the problem is more difficult.

From the point of view of the stability, Gauss-Legendre Runge-Kutta methods are a very good option because they are $A$-stable and $P$-stable, that is, the stability function at pure imaginary points has unit modulus. Therefore, they are stable for any stepsize $h$ and, moreover, they are zero dissipative.  They have high order of accuracy and they are symplectic, an important property when solving Hamiltonian problems \cite{SSC1994,GNI2006}.  A more general class of methods are Hamiltonian Boundary Value Methods (HBVMs) \cite{BIT2010,LIMbook2016}.  They are also $A$-stable, $P$-stable and, moreover, they are energy--conserving.  These two classes of methods have the inconvenient that they are fully implicit and can require a high computational cost. Also, since they are based on Taylor expansions, the local truncation error will depend on $(\omega h)^{p+1}$ and the stepsize could be restricted by accuracy reasons. Nevertheless, HBVMs can be also regarded as spectral methods along the orthonormal Legendre polynomial basis \cite{BIT2012}, and this opens a new perspective in their application. Early references on the usage of spectral methods in time are \cite{BeSt2000,Bo1997,Hu1972,Hu1972-1}, and a further related reference is \cite{TaSu2012}.  In this paper,  we shall use HBVMs as spectral methods also considering a very efficient implementation of such methods, when solving problem (\ref{first}) (or (\ref{sec})).

With these premises, the paper is organized as follows: in Section~\ref{spectral}  we consider the use of the Legendre basis to define a spectral method in time; in Section~\ref{HBVMs}  we see that, by approximating the involved integrals via a Gaussian quadrature, one retrieves HBVMs; in Section~\ref{blendsec} we consider a very efficient nonlinear iteration for solving the generated discrete problems; in Section~\ref{numtest} we present some numerical tests; at last, in Section~\ref{fine} we report a few conclusions and remarks.

\section{Spectral methods}\label{spectral}

We shall here consider, as a suitable orthonormal basis for representing the solution of (\ref{sec2}) (or, more in general, of (\ref{first})) on the interval $[0,h]$, the orthonormal basis for $L^2[0,1]$ functions given by Legendre polynomials:
\begin{equation}\label{Legpol}
P_i\in\Pi_i, \qquad \int_0^1 P_i(x)P_j(x)\dd x = \delta_{ij,}, \qquad \forall i,j=0,1,\dots.
\end{equation}
We start considering the approximation of the linear part in (\ref{sec2}). For this purpose, we need the following preliminary results.

\begin{lem}\label{lem1} Let $\omega h>0$. For $s=1,2,\ldots$, set 
\begin{equation}\label{intPs}
g(s,\omega h) ~:=~ \sqrt{\frac{(2s+1)\pi}{\omega h}} \left|J_{s+\frac{1}2}\left( \frac{\omega h}2\right)\right|,
\end{equation}
where $J_{s+\frac{1}2}(\cdot)$ is the Bessel function of the first kind. Then, for all $s=1,2,\ldots:$
\begin{equation}\label{intPs1}
\left|\int_0^1 P_s(c)\cos(\omega h c)\dd c\right|,\, \left|\int_0^1 P_s(c)\sin(\omega h c)\dd c\right| ~\le~g(s,\omega h).
\end{equation}
\end{lem}
\proof The proof easily derives from \cite[Eq. (9)]{EW1999}, which states that
\begin{equation}\label{uguale}
\left|\int_0^1 P_s(c)\cos(\omega h c)\dd c\right|^2+\left|\int_0^1 P_s(c)\sin(\omega h c)\dd c\right|^2 ~=~  g(s,\omega h)^2.
\QED
\end{equation}

\begin{figure}[ht]
\centerline{\includegraphics[width=12cm,height=8cm]{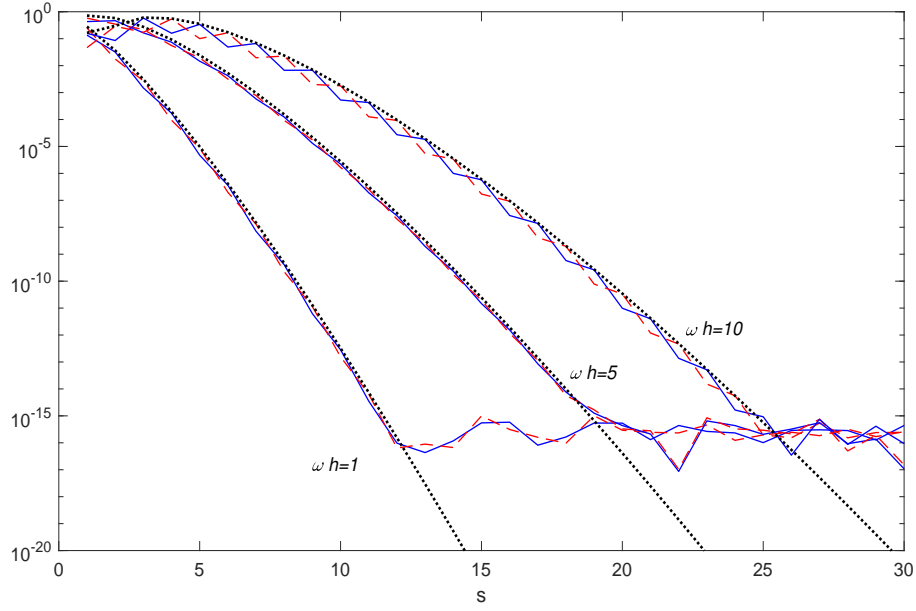}}
\caption{Values of the integrals at the left-hand side in (\ref{intPs1}), numerically evaluated by a high order Gauss-Legendre formula (solid and dashed lines, respectively), along with their bound $g(s,\omega h)$ defined in (\ref{intPs}) (dotted lines), for $\omega h = 1,5,10$.}
\label{essefig}
\end{figure}

\medskip
In Figure~\ref{essefig} we plot the values of the integrals in (\ref{intPs1}), numerically computed via a high-order Gauss-Legendre quadrature formula (solid lines and dashed lines, respectively, for the two integrals at the left-hand side in (\ref{intPs1})), together with the bound provided by the function defined in (\ref{intPs}) (dotted lines), for $\omega h = 1,5,10$. As one may see, they are in very good agreement, until round-off error level is reached, so that the numerical quadrature becomes ineffective.

\begin{lem}\label{lem2} For $\omega h>0$, and $s\gg1$, the function $g(s,\omega h)$ defined in (\ref{intPs}) is an increasing function of $\omega h$ and a decreasing function of $s$.
\end{lem}
\proof From \cite[Eq. (20)]{EW1999}, one derives that, for $\omega h>0$ and $s\gg 1$,
\begin{equation}\label{larges}
g(s,\omega h) ~=~ \sqrt{\frac{(2s+1)\pi}{\omega h}} \left|J_{s+\frac{1}2}\left( \frac{\omega h}2\right)\right| ~\approx~ \sqrt{\frac{e}{2(2s+1)}}\left(\frac{e\omega h}{2(2s+1)}\right)^s.
\end{equation}
The latter function, in turn, is an increasing function of $\omega h$, for fixed $s$, and a decreasing function of $s$, for any fixed $\omega h>0$ and  all $s$ such that $2(2s+1)> e\omega h$.\QED

\medskip
We observe that also the result of Lemma~\ref{lem2} is clearly confirmed by the plots in Figure~\ref{essefig}.

\medskip
Let now consider the approximation of the linear part of the problem (\ref{sec2}) on the interval $[0,h]$, i.e.
\begin{equation}\label{firstlin}
\dot y = J_2\otimes A\, y, \qquad y(0)=y_0:=\pmatrix{c}q_0\\ p_0\endpmatrix.
\end{equation}
The following result holds true.

\begin{theo}\label{exp} The solution of (\ref{firstlin}) satisfies
\begin{equation}\label{expy}
y(ch) ~\equiv~ \ee^{J_2\otimes A hc}y_0 ~=~ y_0 + h\sum_{j\ge 0} \int_0^c P_j(x)\dd x\, \gamma_j (y), \qquad c\in [0,1],
\end{equation}
with
\begin{equation}\label{gammaj}
\gamma_j(y) =  J_2\otimes A\int_0^1 P_j(\tau) y(\tau h)\dd \tau,\qquad j=0,1,2,\dots.  
\end{equation}
Moreover, 
\begin{equation}\label{y1y0}
y(h)=y_0 + h\gamma_0(y).
\end{equation}
\end{theo}
\proof According to \cite{BIT2012}, the right-hand side of (\ref{firstlin}) can be expanded, on the interval $[0,h]$, along the Legendre basis:
$$\dot y(ch) =  \sum_{j\ge0} P_j(c) \gamma_j(y), \qquad c\in[0,1],$$
with the coefficients $\gamma_j(y)$ clearly given by (\ref{gammaj}), because of the orthonormality conditions (\ref{Legpol}). Integration side by side of such equation, and imposing the initial condition, then gives (\ref{expy}). At last, (\ref{y1y0}) follows from (\ref{expy}), by considering that $$\int_0^1 P_j(x)\dd x = \delta_{j0}.\,\QED$$

\medskip
By further considering that 
$$\ee^{J_2\otimes A} = \pmatrix{cc} \cos(A) &\sin(A)\\ -\sin(A) &\cos(A)\endpmatrix,$$
one has then that $y$ in (\ref{expy}) is also given by  
\begin{equation}\label{ych}
y(ch) = \pmatrix{cc} \cos(A hc) &\sin(A hc)\\ -\sin(Ahc) &\cos(Ahc)\endpmatrix y_0, \qquad c\in[0,1],
\end{equation}
namely, $y$ is obtained as the combination of sines and cosines, with frequencies not larger than
\begin{equation}\label{omega}
\omega = \|A\|.
\end{equation}
As a consequence of Theorem~\ref{exp} and (\ref{omega}), one has the following.

\begin{crite}\label{crits0} 
When using a finite precision arithmetic with machine epsilon $u$,  and with reference to the function $g$ defined in (\ref{intPs}), the series at the rigt-hand side in (\ref{expy}) can be truncated at a convenient value $s_0-1$ such that:
\begin{equation}\label{esse}
{\blue g(s_0,\omega h) < u \cdot \max_{j< s_0} \, g(j,\omega h).}
\end{equation}
\end{crite}

\smallskip
The derivation of such criterion is as follows. From the results of Lemmas~\ref{lem1} and \ref{lem2},  (\ref{Asdp}), (\ref{gammaj}), (\ref{ych}), and with reference to the function $g$ defined in (\ref{intPs}),
\begin{eqnarray*} 
\gamma_j(y) &=& J_2\otimes A\int_0^1 P_j(c) \pmatrix{cc} \cos(A hc) &\sin(A hc)\\ -\sin(Ahc) &\cos(Ahc)\endpmatrix y_0\,\dd c \\
&=& (J_2\otimes Q)  \underbrace{\int_0^1 P_j(c) \pmatrix{cc} \Lambda\cos(\Lambda hc) &\Lambda\sin(\Lambda hc)\\ -\Lambda\sin(\Lambda hc) &\Lambda\cos(\Lambda hc)\endpmatrix \,\dd c}_{=:M_j}\, (I_2\otimes Q^\top) y_0\\
&=:& (J_2\otimes Q) \, M_j \, (I_2\otimes Q^\top) y_0,
\end{eqnarray*}
where $I_2$ is the $2\times 2$ identity matrix. By considering the 2-norm, we shall then consider the approximation
\begin{equation}\label{appr1}
\|\gamma_j(y)\|_2 \sim \|M_j\|_2 \|y_0\|_2.
\end{equation}
Moreover, from (\ref{uguale}), one has:
\begin{eqnarray}\nonumber
\|M_j\|_2^2 &=& \max_{\lambda\in\sigma(A)}  \lambda^2 \left[\Big(\int_0^1 P_j(\tau)\cos(\lambda h\tau)\dd\tau\Big)^2+
\Big(\int_0^1 P_j(\tau)\sin(\lambda h\tau)\dd\tau\Big)^2\right]\\ \label{appr2}
&=& \max_{\lambda\in\sigma(A)}  \left[ \lambda g(j,\lambda h) \right] ^2 ~{\blue \le}~ \left[ \omega g(j,\omega h)\right]^2.
\end{eqnarray}
{\blue Clearly, the last inequality in (\ref{appr2}) becomes an equality, if in (\ref{omega}) one considers the $\|\cdot\|_2$.}
By taking into account that, according to Lemma~\ref{lem2}, for $j\gg1$ the function $g(j,\omega h)$ is a decreasing function of $j$, and considering that
$$\left|\int_0^c P_j(x)\dd x\right| \le \frac{1}{\sqrt{2j+1}}, \qquad c\in[0,1], \qquad j\ge0,$$
one can neglect the terms in the series (\ref{expy}), starting from the index $s_0$ such that
$$\|\gamma_{s_0}(y)\|_2 < u\cdot \max_{j<s_0} \|\gamma_j(y)\|_2.$$
The criterion (\ref{esse}) then follows from (\ref{appr1}) and (\ref{appr2}), by considering the estimate
$$\|\gamma_j(y)\|_2 \sim \omega\, g(j,\omega h) \|y_0\|_2.\,\QED$$\medskip

By means of Criterion~\ref{crits0}, from (\ref{expy}) one obtains
\begin{equation}\label{ys}
y(ch) ~\doteq~ y_{s_0}(ch) ~:=~y_0 + h\sum_{j=0}^{s_0-1} \int_0^c P_j(x)\dd x\, \gamma_j(y_{s_0}), \qquad c\in [0,1],
\end{equation}
where $\gamma_j(y_{s_0})$ is formally still defined by (\ref{gammaj}), by replacing $y$ with $y_{s_0}$. In the above expression and hereafter, \,$\doteq$\, means ``equal within round-off error level''. In fact, neither the terms starting from $s_0$  would be taken into account by the used finite precision arithmetic, nor they could be numerically reliably computed, as is confirmed by the plots in Figure~\ref{essefig}.

\begin{figure}[t]
\centerline{\includegraphics[width=12cm,height=8cm]{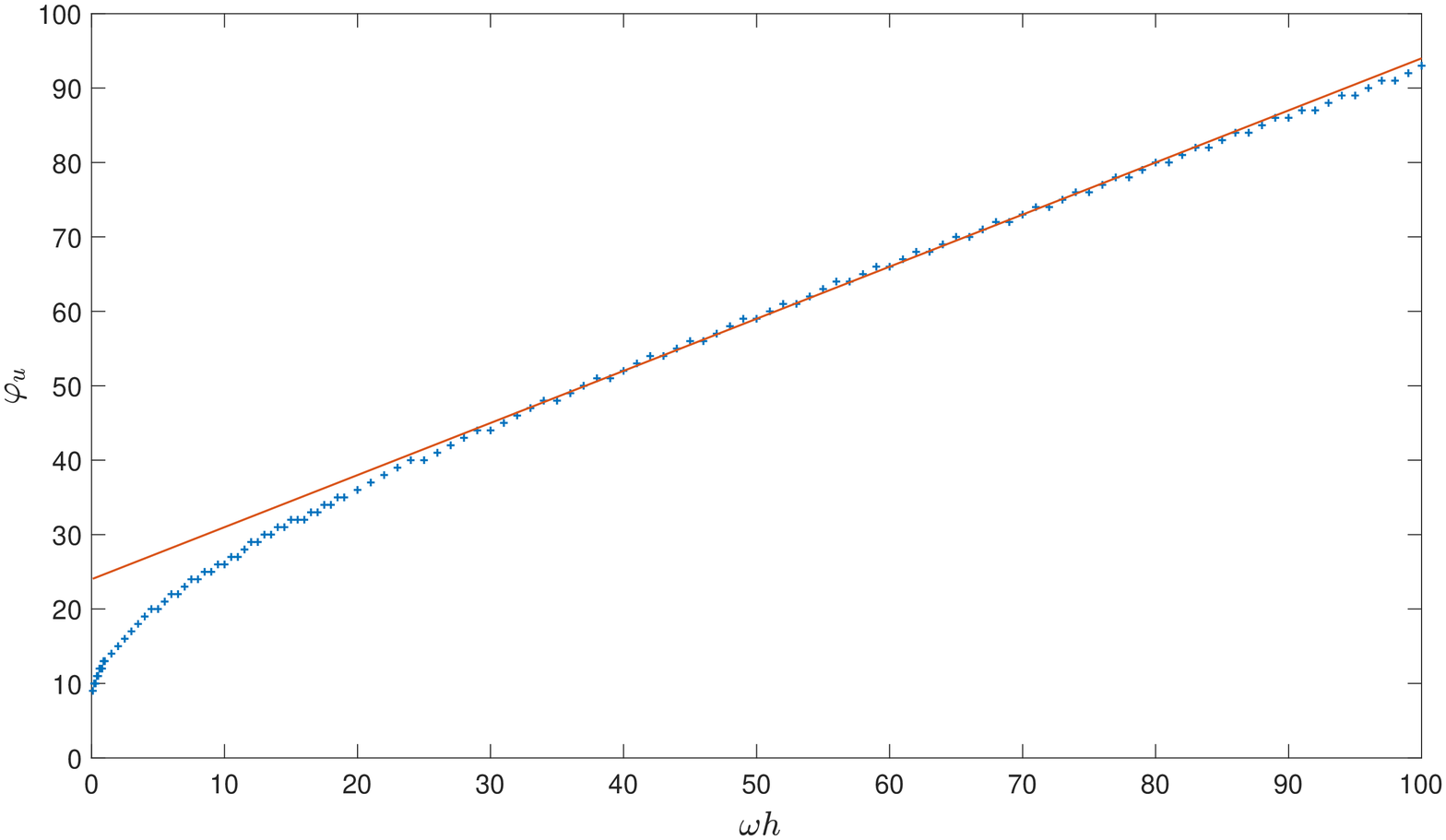}}
\caption{function $\varphi_u$ versus $\omega h$, for the double precision IEEE. For $\omega h\gg1$ one has that $\varphi_u(\omega h)\approx 24+0.7\cdot\omega h$.}
\label{esse0fig}
\end{figure}

In addition to this, from Criterion~\ref{crits0}, one clearly obtains that $s_0$ is provided by a function, say $\varphi_u$, of $\omega h$ and also depending on the used machine epsilon $u$:
\begin{equation}\label{varfi}
s_0 = \varphi_u(\omega h).
\end{equation} 
In Table~\ref{stab} we list a few values of  $s_0$, depending on the product $\omega h$, when considering the double precision IEEE.
Moreover, in Figure~\ref{esse0fig} we plot $\varphi_u(\omega h)$ versus $\omega h$, for the double precision IEEE. From the figure, one has that
$$\varphi_u(\omega h) \approx 24+0.7\cdot \omega h, \qquad \omega h\gg 1.$$

\begin{table}[t]
\caption{Values of the parameter $s_0$ in (\ref{esse}), as a function of $\omega h$, for the double precision IEEE.}

\smallskip
\label{stab}
\centerline{
\begin{tabular}{|r|ccccccccc|}
\hline
$\omega h$  & 0.1   & 0.5    & 1       &  5 & 10     &   25 &  50    &  75 &  100\\
\hline
$s_0$           & 9     & 11 &    13 & 20 & 26 & 40 & 59 & 76 & 93 \\
\hline
\end{tabular}}
\end{table}

Next, let us consider the complete problem (\ref{sec1})--(\ref{sec2}). By assuming the ansatz
\begin{equation}\label{ansatz}
\nabla f(q(t)) ~\sim~  \cos(\nu At) \tilde q_0, \qquad t\in[0,h],
\end{equation}
for suitables $\tilde q_0\in\RR^m$ and $\nu>1$.\footnote{I.e., locally $\nabla f$ approximately behaves as a polynomial of degree $\nu$.}
Consequently, by the nonlinear variation of constants formula, one obtains that the solution of (\ref{sec2}), on the interval $[0,h]$, can be approximately expressed as
\begin{eqnarray}\label{vcf}
q(ch) &\sim& \cos( A\/ch )q_0 + \sin(A\/ch)p_0 - h A^{-1}\int_0^c \sin (Ah\/(c-\tau)) \cos(\nu Ah\tau)\dd\tau\, \tilde q_0,\\ \nonumber
p(ch) &\sim& -\sin( A\/ch )q_0 + \cos(A\/ch)p_0- h A^{-1}\int_0^c \cos (Ah\/(c-\tau)) \cos(\nu Ah\tau)\dd\tau\, \tilde q_0, \quad c\in[0,1].
\end{eqnarray}
Further, by considering that, by the Werner formulae,
\begin{eqnarray*}
\sin (Ah\/(c-\tau)) \cos(\nu Ah\tau) &=& \frac{1}2\left[ \sin\left( Ah (c+(\nu-1)\tau)\right) + \sin\left( Ah (c-(\nu+1)\tau)\right)\right],\\
\cos (Ah\/(c-\tau)) \cos(\nu Ah\tau)&=& \frac{1}2\left[ \cos\left( Ah (c+(\nu-1)\tau)\right) + \cos\left( Ah (c-(\nu+1)\tau)\right)           \right],
\end{eqnarray*}
so that
\begin{eqnarray*}
\lefteqn{
\int_0^c \sin (Ah\/(c-\tau)) \cos(\nu Ah\tau)\dd\tau}\\ &=&\frac{1}2\int_0^c\left[ \sin\left(  Ah (c+(\nu-1)\tau)\right) + \sin\left( Ah (c-(\nu+1)\tau)\right)\right]\dd\tau\\
&=&A^{-1}h^{-1} \frac1{\nu^2-1}\left[ \cos\left(hAc\right) - \cos\left(hA\nu c\right)\right] ,
\end{eqnarray*}
and, similarly,
\begin{eqnarray*}
\lefteqn{
\int_0^c \cos (Ah\/(c-\tau)) \cos(\nu Ah\tau)\dd\tau}\\ &=&\frac{1}2\int_0^c\left[ \cos\left(  Ah (c+(\nu-1)\tau)\right) + \cos\left(  Ah (c-(\nu+1)\tau)\right)\right]\dd\tau\\
&=&A^{-1}h^{-1} \frac1{\nu^2-1}\left[ \nu\sin\left(hA\nu c\right) - \sin\left(hAc\right)\right] ,
\end{eqnarray*}
the following result is proved.

\begin{theo} Under the ansatz (\ref{ansatz}) the solution (\ref{vcf}) of (\ref{sec1})--(\ref{sec2}) has an oscillatory behavior, with component in the maximum frequency proportional to:
\begin{equation}\label{wmax1}
\frac{A^{-2}}{\nu^2-1}\left[\cos\left(hA\nu c\right) -\nu\sin\left(hA\nu c\right) \right]\tilde q_0, \qquad c\in[0,1].
\end{equation}
\end{theo}

As a consequence of the previous result, by expanding again $y(t)=(q(t),p(t))^\top$ in (\ref{sec2}) along the Legendre basis,
\begin{equation}\label{expy1}
y(ch) ~=~ y_0 + h\sum_{j\ge 0} \int_0^c P_j(x)\dd x  \,\psi_j(y), \qquad c\in [0,1],
\end{equation}
with  
\begin{equation}\label{psij}
\psi_j(y) =  \int_0^1 P_j(\tau)\left[J_2\otimes A \,y(\tau h)) + J_2\otimes I_m {\blue \tilde{f}(y(\tau h)})\right]\dd\tau,\qquad j=0,1,2,\dots,
\end{equation}
one derives the following criterion, which generalises the one given by Criterion\,\ref{crits0}, for the linear problem (\ref{firstlin}), to the whole nonlinear problem (\ref{sec1})--(\ref{sec2}).

\begin{crite}\label{crits}
When using a finite precision arithmetic with machine epsilon $u$,  and with reference to the function $g$ defined in (\ref{intPs}), the series at the rigt-hand side in (\ref{expy1}) can be truncated at a convenient value $s-1$ such that:
\begin{equation}\label{esse1}
{\blue g(s,\nu \omega h) < u \cdot \max_{j< s} \, g(j,\nu \omega h).}
\end{equation}
\end{crite}
Consequently, one obtains full machine accuracy by considering, in place of (\ref{expy1}):
\begin{equation}\label{expy1s}
y(ch) ~\doteq~ y_s(ch) ~:=~y_0 + h\sum_{j=0}^{s-1} \int_0^c P_j(x)\dd x  \,\psi_j(y_s), \qquad c\in [0,1],
\end{equation}
where the coefficients $\psi_j(y_s)$ are formally defined as in (\ref{psij}), by replacing $y$ by $y_s$.  We observe that, similarly to what seen for the parameter $s_0$, now 
$$s = \varphi_u(\nu\omega h),$$
where $\varphi_u$ is the same function defined in (\ref{varfi}).\footnote{Clearly, when $\nu=1$, then $s$ coincides with $s_0$, as defined in (\ref{esse}).}

Finally, we observe that, when solving problem (\ref{first}), instead of (\ref{sec2}), one would formally obtain the same relation (\ref{expy1s}), with the coefficients $\psi_j(y)$ defined as:
\begin{equation}\label{psijuno}
\psi_j(y) =  \int_0^1 P_j(\tau) J\left[ A y(\tau h) + \nabla f(y(\tau h))\right]\dd\tau,\qquad j=0,1,2,\dots,
\end{equation}
in place of (\ref{psij}). Also in such a case, the Criteria~\ref{crits0} and \ref{crits} continue formally to hold.

In the next section, we study the properties of the numerical method obtained by approximating to full machine accuracy the integrals appearing in (\ref{psij}) or (\ref{psijuno}) by means of a suitable quadrature rule. In particular, we shall consider a Gaussian quadrature based at the $k$ Legendre abscissae (thus, of order $2k$) where, in order to guarantee full machine accuracy, when using the IEEE double precision, we choose
\begin{equation}\label{kappa}
k = \max\{ s+2, 20\},    
\end{equation}
with $s$ defined according to (\ref{esse1}). However, more refined choices could in principle be used.

\section{Hamiltonian Boundary Value Methods (HBVMs)}\label{HBVMs} 
In this section we illustrate the application of Hamiltonian Boundary Value Methods (HBVMs) for solving problem (\ref{first}), which is more general than (\ref{sec2}). HBVMs form a class of energy-conserving Runge-Kutta methods which has been studied in a series of papers  for the numerical solution of Hamiltonian problems \cite{BIT2010,BIT2011,BIT2012_0,BIT2012,BFCI2014,BIT2015}. They have been also extended along a number of directions   \cite{BIT2012_1,BIT2012_2,BrIa2012,BCMR2012,BFCI2015,ABI2015,BBFCI2018} (we also refer to the recent monograph \cite{LIMbook2016}).

In more details, a HBVM$(k,s)$ method is the $k$-stage Runge-Kutta method defined by the Butcher tableau
\begin{equation}\label{butab}
\begin{array}{c|c} \bfc & \I_s \P_s^\top \Omega \\ \hline \\[-3mm] & \bfb^\top\end{array}, \qquad \bfb = \pmatrix{ccc} b_1 &\dots & b_k\endpmatrix ^\top, \quad \bfc = \pmatrix{ccc} c_1 &\dots & c_k\endpmatrix ^\top,
\end{equation}
where, with reference to the Legendre polynomial basis defined in (\ref{Legpol}), $(c_i,b_i)$ are the abscissae and weights of the Gauss-Legendre quadrature formula of order $2k$ (i.e., $P_k(c_i)=0$, $i=1,\dots,k$),
and
\begin{eqnarray}\nonumber
\P_s &=& \pmatrix{ccc} 
P_0(c_1) & \dots &P_{s-1}(c_1)\\
\vdots      &          &\vdots\\
P_0(c_k) & \dots & P_{s-1}(c_k)
\endpmatrix \in\RR^{k\times s}, \\[2mm] \label{PIO}
 \I_s &=& \pmatrix{ccc} 
\int_0^{c_1} P_0(x)\dd x & \dots &\int_0^{c_1} P_{s-1}(x)\dd x\\
\vdots      &          &\vdots\\
\int_0^{c_k} P_0(x)\dd x & \dots & \int_0^{c_k} P_{s-1}(x)\dd x
\endpmatrix \in\RR^{k\times s},\\[2mm] \nonumber
 \Omega &=& \pmatrix{ccc} b_1 \\ &\ddots \\ && b_k\endpmatrix \in\RR^{k\times k}.
 \end{eqnarray}
The following result is known to hold for such methods (see, e.g., \cite{BIT2012,LIMbook2016}).

\begin{theo}\label{HBVMth} For all $k\ge s$, the HBVM$(k,s)$ method (\ref{butab}):
\begin{itemize}
\item is symmetric and has order $2s$;
\item when $k=s$ it reduces to the (symplectic) $s$-stage Gauss collocation method;
\item it is energy-conserving, when applied for solving Hamiltonian problems with a polynomial Hamiltonian of degree not larger than
$2k/s$;
\item for general and suitably regular Hamiltonians, the Hamiltonian error per step is $O(h^{2k+1})$.
\end{itemize}
\end{theo} 

\begin{rem} Because of the result of Theorem~\ref{HBVMth}, one has that an exact energy-conservation is always obtained, by choosing $k$ large enough, in the polynomial case. Moreover, even in the non-polynomial case, a {\em practical } energy-conservation can always be gained, by choosing $k$ large enough so that the Hamiltonian error falls within the round-off error level. This, in turn, doesn't affect too much the computational cost of the method, as is shown in Section~\ref{blendsec} (see also \cite{BIT2011,BFCI2014,LIMbook2016}). In particular, the choice (\ref{kappa}) of $k$ will always provide us with a practical energy conservation.
\end{rem}

For sake of completeness, we also mention that one may consider the limit as $k\rightarrow\infty$ of HBVMs \cite{BIT2010}, thus obtaining a continuous-stage Runge-Kutta method (see also \cite{BIT2010,Ha2010} and \cite[Chapter\,3.5]{LIMbook2016}).\medskip

A few properties of the matrices defined in (\ref{PIO}) are here recalled, for later use.

\begin{lem}\label{PIOX} For all $k\ge s$, one has:
\begin{equation}\label{Xs}
\P_s^\top\Omega \I_s = X_s \equiv \pmatrix{cccc}
\xi_0 & -\xi_1\\
\xi_1 & 0  & \ddots \\
         &\ddots &\ddots &-\xi_{s-1}\\
         &           &\xi_{s-1} &0\endpmatrix, \qquad \xi_i=\frac{1}{2\sqrt{|4i^2-1|}}, \quad i=0,1,\dots.
\end{equation}         
 Moreover,
$$\det(X_s) = \left\{ \begin{array}{ccl}
\prod_{i=1}^{s/2} \xi_{2i-1}^2, & ~&s\mbox{~even},\\ ~\\
\prod_{i=0}^{(s-1)/2} \xi_{2i}^2, & ~&s\mbox{~odd},
\end{array}\right.$$
so that matrix $X_s$ is nonsingular, for all $s=1,2,\dots$.
\end{lem}
\proof See, e.g., \cite[Lemmas~3.6 and 3.7]{LIMbook2016}.\,\QED
\medskip

Let now study the application of a HBVM$(k,s)$ method, $k>s$, for solving (\ref{first}). In so doing,  
by setting $$Y :=\pmatrix{c} Y_1\\ \vdots \\ Y_k\endpmatrix,  \qquad \nabla f(Y) := \pmatrix{c}\nabla f( Y_1)\\ \vdots \\ \nabla f(Y_k)\endpmatrix, \qquad \bfe = \pmatrix{c}1\\ \vdots\\1\endpmatrix\in\RR^k,$$
with $Y$ the stage vector of the method, one obtains the discrete problem of block dimension $k$,
\begin{equation}\label{Yp}
Y = \bfe\otimes y_0 + h \I_{s}\P_{s}^\top\Omega \otimes J
\left[  I_k\otimes A \,Y + \nabla f(Y)\right], 
\end{equation}
and the following approximation to $y(h)$:
\begin{equation}\label{y1new}
y_1 = y_0 + h\sum_{i=1}^k b_i {\blue J[A Y_i + \nabla f(Y_i)].}
\end{equation}
Nevertheless, by considering that
\begin{equation}\label{Yi}
Y_i = y_0 + h\sum_{j=0}^{s-1} \int_0^{c_i} P_j(x)\dd x\, \psi_j ~=:~\sigma_{s}(c_ih), \qquad i=1,\dots,k,
\end{equation}
with $\sigma_{s}\in\Pi_{s}$, and
\begin{eqnarray}\nonumber
\psi_j &=& \sum_{\ell=1}^k b_\ell P_j(c_\ell) J[A Y_\ell + \nabla f(Y_\ell)]\\
\label{psij1}
&\equiv&\sum_{\ell=1}^k b_\ell P_j(c_\ell) J[A \sigma_{s}(c_\ell h) + \nabla f(\sigma_{s}(c_\ell h))], \qquad j=0,\dots,s-1,
\end{eqnarray}
comparison of (\ref{Yi})--(\ref{psij1}) with  (\ref{expy1s}), provides us with the following result.

\begin{theo}\label{klarge} Provided that the quadrature is exact within full machine accuracy, with reference to (\ref{expy1s}) and (\ref{psijuno}), one has:
$$\sigma_{s}(ch) \doteq y_{s}(ch) \quad\Rightarrow\quad \psi_j \doteq \psi_j(y_{s}), \quad j=0,\dots,s-1.$$
\end{theo}

\begin{rem}\label{norder}
The result of Theorem~\ref{klarge}, with $s$ chosen according to Criterion~\ref{crits}, clearly shows that for such a method the concept of order doesn't apply, since it will always provide the maximum possible accuracy, for the used finite precision arithmetic.
\end{rem}

In order to improve the computational efficiency, instead of directly solving the discrete problem (\ref{Yp}), having block dimension $k$, we shall consider a more convenient formulation of it \cite{BIT2011}. In more details, by setting
\begin{equation}\label{bfpsi}
\bfpsi := \pmatrix{c} \psi_0\\ \vdots \\ \psi_{s-1}\endpmatrix \equiv \P_{s}^\top \Omega \otimes J \left[  (I_k\otimes A) \,Y + \nabla f(Y)\right],
\end{equation}
from (\ref{Yp}) one obtains, 
$$Y = \bfe\otimes y_0 + h\I_{s}\otimes I_{2m} \bfpsi,$$
which, plugged into (\ref{bfpsi}), provides the following new discrete problem, clearly equivalent to (\ref{Yp}):
\begin{equation}\label{bfpsip}
G(\bfpsi) \,:=\, \bfpsi - \P_{s}^\top\Omega \otimes J\left[ (I_k\otimes A) \left(\bfe\otimes y_0 + h\I_{s}\otimes I_{2m} \bfpsi \right) +\nabla f \left(\bfe\otimes y_0 + h\I_{s}\otimes I_{2m} \bfpsi \right) \right]\,=\,\bf0.
\end{equation}
Once it has been solved, the new approximation (\ref{y1new}) turns out to be given by:
\begin{equation}\label{y1p}
y_1 = y_0 + h\psi_0 \equiv \sigma_s(h) \doteq y_{s}(h).
\end{equation}
We observe that the advantage of solving the discrete problem (\ref{bfpsip}) instead of the stage problem (\ref{Yp}) is twofold:
\begin{enumerate}
\item the problem (\ref{bfpsip}) has block dimension $s$, independently of $k$;
\item its numerical solution via a simplified Newton-type iteration is very efficient.
\end{enumerate}
The last point is elucidated in the next section.

\section{Efficient implementation of the methods}\label{blendsec}

The use of the simplified Newton method for solving (\ref{bfpsip}) is described, by virtue of (\ref{Xs}), by the following iteration:
\begin{eqnarray}\nonumber
\mbox{initialize}&& \bfpsi^0\\ \label{newt0}
\mbox{for~}\ell&=&0,1,\ldots: \\ \nonumber
&\mbox{solve} &[I-hX_{s}\otimes J(A+\nabla^2 f(y_0))] \Delta^\ell = -G(\bfpsi^\ell)\\ \nonumber
&\mbox{set}& \bfpsi^{\ell+1} = \bfpsi^\ell + \Delta^\ell\\ \nonumber
\mbox{end}
\end{eqnarray}
with $I$ the identity matrix of dimension $s\cdot 2m$.

\smallskip
Firstly, we notice that, by virtue of (\ref{a2large}), we can consider the approximation
\begin{equation}\label{bastaA}
 A+\nabla^2 f(y_0) \approx A,
\end{equation} 
thus obtaining a coefficient matrix which is {\em constant} and is the same for all the integration steps.
Consequently, the iteration (\ref{newt0}) simplifies to:
\begin{eqnarray}\nonumber
\mbox{initialize}&& \bfpsi^0\\  \label{newt1}
\mbox{for~}\ell&=&0,1,\ldots: \\ \nonumber
&\mbox{solve} &[I-hX_{s}\otimes JA ]\, \Delta^\ell = -G(\bfpsi^\ell)\\ \nonumber
&\mbox{set}& \bfpsi^{\ell+1} = \bfpsi^\ell + \Delta^\ell\\ \nonumber
\mbox{end}
\end{eqnarray}
This iteration can be further simplified by using a Newton-splitting {\em blended iteration}. This technique, at first devised in \cite{Br2000,BrMa2002}, has then be generalized \cite{BrMa2007,BrMa2009} and  implemented in the computational codes {\tt BiM} \cite{BrMa2004} and {\tt BiMD} \cite{BrMaMu2006} for stiff ODE/DAE IVPs. It has been also considered for HBVMs \cite{BIT2011} and is implemented in the Matlab code {\tt HBVM} \cite{LIMbook2016}. The novelty, in the present case, is due to the approximation (\ref{bastaA}), which makes it extremely efficient. As a result, the iteration (\ref{newt1}) modifies as follows:
\begin{eqnarray}\nonumber
\mbox{initialize}&& \bfpsi^0\\ \label{blend} 
\mbox{for~}\ell&=&0,1,\ldots: \\ \nonumber
&\mbox{set}     & \bfeta^\ell = -G(\bfpsi^\ell)\\ \nonumber
&\mbox{set}     & \bfeta_1^\ell = \left[(\rho_s X_s^{-1})\otimes I_{2m}\right] \bfeta^\ell \\ \nonumber
&\mbox{set}     & \bfu^\ell = \left[I_s\otimes \Sigma\right] (\bfeta^\ell -\bfeta_1^\ell)\\ \nonumber 
&\mbox{set}     & \Delta^\ell = \left[I_s\otimes \Sigma\right] (\bfeta_1^\ell +\bfu^\ell)\\ \nonumber
&\mbox{set}     & \bfpsi^{\ell+1} = \bfpsi^\ell + \Delta^\ell  \\ \nonumber
\mbox{end}
\end{eqnarray}
Here, according to \cite{BrMa2002,BrMa2009,BIT2011},
\begin{equation}\label{ros}
\rho_s = \min_{\lambda\in\sigma(X_s)} |\lambda|,
\end{equation}
and 
\begin{equation}\label{sigma}
\blue \Sigma = \left( I_{2m} - h\rho_s JA \right)^{-1}.
\end{equation}
Consequently, we notice that one needs only to compute {\em once and for all} matrix $\Sigma$ (or factorize $\Sigma^{-1}$), having the same size as that of the continuous problem. 

\smallskip
Secondly, in order to gain convergence for relatively large values of $\omega h$, it is important to choose an appropriate starting value $\bfpsi^0$ in (\ref{blend}). For this purpose, we use the solution of the associated homogeneous problem (\ref{firstlin}).
This latter can be conveniently computed, by virtue of Criterion~\ref{crits0} (see (\ref{esse})) by the HBVM$(s_0,s_0)$ method (i.e., the $s_0$-stage Gauss collocation method). Consequently, by repeating similar steps as above, by defining the vectors
$$\bfgamma = \pmatrix{c} \gamma_0\\ \vdots \\ \gamma_{s_0-1}\endpmatrix \in\RR^{s_02m},\qquad 
\tilde\bfe = \pmatrix{c} 1\\ \vdots \\ 1\endpmatrix\in \RR^{s_0},$$ 
and the matrices, with a structure similar to (\ref{PIO}),
$$\tilde\P_{s_0} = \left( P_{j-1}(\tilde{c}_i) \right), ~\tilde\I_{s_0} = \left( \int_0^{\tilde{c}_i} P_{j-1}(x)\dd x \right), ~ \tilde\Omega= \pmatrix{ccc} \tilde{b}_1 \\ &\ddots \\&& \tilde{b}_{s_0}\endpmatrix~\in\RR^{s_0\times s_0},$$ with $(\tilde{c}_i,\tilde{b}_i)$ the abscissae and weights of the Gauss-Legendre quadrature of order $2s_0$,
one solves the discrete problem
$$\tilde{G}(\bfgamma) := \bfgamma -\tilde\P_{s_0}^\top \tilde\Omega \otimes JA \left[ \tilde\bfe\otimes y_0 + h\tilde\I_{s_0}\otimes I_{2m} \bfgamma\right] = \bf0.$$
This, in turn, can be done by means of the following approximate {\em blended} iteration, similar to (\ref{blend}):
\begin{eqnarray}\nonumber
\mbox{initialize}&& \bfgamma^0 = \bf0\\  \label{blend0}
\mbox{for~}\ell&=&0,1,\ldots: \\ \nonumber
&\mbox{set}     & \tilde\bfeta^\ell = -\tilde{G}(\bfgamma^\ell)\\ \nonumber
&\mbox{set}     & \tilde\bfeta_1^\ell = \left[(\rho_{s_0} X_{s_0}^{-1})\otimes I_{2m}\right] \tilde\bfeta^\ell \\ \nonumber
&\mbox{set}     & \tilde\bfu^\ell = \left[I_{s_0}\otimes \Sigma\right] (\tilde\bfeta^\ell -\tilde\bfeta_1^\ell)\\ \nonumber 
&\mbox{set}     & \tilde\Delta^\ell = \left[I_{s_0}\otimes \Sigma\right] (\tilde\bfeta_1^\ell +\tilde\bfu^\ell)\\ \nonumber
&\mbox{set}     & \bfgamma^{\ell+1} = \bfgamma^\ell + \tilde\Delta^\ell\\ \nonumber
\mbox{end}
\end{eqnarray}
where matrix $X_{s_0}\in\RR^{s_0\times s_0}$ is defined according to (\ref{Xs}), and the same matxix $\Sigma$ defined in (\ref{sigma}) can be used.\footnote{The approximation stems from the fact that in (\ref{sigma}) $\rho_s$ is used in place of $\rho_{s_0}$.} Once this has been done, the following initialization can be conveniently used for the iteration (\ref{blend}):
\begin{equation}\label{bfpsi0}
\bfpsi^0 = \pmatrix{c} \bfgamma \\ {\bf0}\endpmatrix, \qquad\mbox{with}\qquad  {\bf0}\in\RR^{(s-s_0)2m}.  
\end{equation}
\begin{defi}
We shall refer to the method defined by the iterations (\ref{blend}) and (\ref{blend0})--(\ref{bfpsi0}) to as {\em spectral HBVM} with parameters $k,s,s_0$, in short SHBVM$(k,s,s_0)$.
\end{defi}

\section{Numerical tests}\label{numtest}

In this section we consider a few numerical tests, aimed at assessing the effectiveness of the SHBVM$(k,s,s_0)$ method, implementing the spectral-Legendre methods described in Section~\ref{spectral}. This will be done by comparing such methods with known existing ones. In particular, we consider the methods below specified, where the first three methods can be used only for problems in the form  (\ref{sec}).

\subsubsection*{Methods:}
\begin{description}
\item[St\"ormer-Verlet]: the St\"ormer-Verlet method \cite{GNI2006}, of order 2;
\item[Gautschi]: the trigonometric fitted method of Gautschi \cite{Ga1961}, of order 2;
\item[Deuflhard]: the improved trigonometric fitted method of Deuflhard \cite{De1979}, of order 2; 
\item[Expode]: the exponential integrator code EXPODE  \cite{Jan2014}, which is based on exponential integrators.  We have used the routine {\tt exprk} based on an exponential RK methods of order 4 described in \cite{HoOs2005};
\item[Gauss]: the $s$-stage Gauss method, $s=1,\dots,4$, of order $2s$;
\item[SHBVM]: the spectral method SHBVM$(k,s,s_0)$.
\end{description}
All numerical tests have been performed on a laptop with {\blue a 2.2\,GHz} dual-core $i7$ and 8GB of memory, running Matlab (R2017b).

\subsubsection*{Duffing equation}
The first test problem is the Duffing equation:
\begin{equation}\label{dufeq}
\ddot q = -(\kappa^2+{\blue\beta^2}) q +2\kappa^2 q^3, \quad t>0, \qquad q(0)=0, \quad \dot q(0) = \beta,
\end{equation}
with Hamiltonian 
\begin{equation}\label{dufH}
H(q,p) = \frac{1}2\left[ p^2 + (\kappa^2+\beta^2) q -\kappa^2 q^4\right], \qquad p=\dot q.
\end{equation} 
In such a case, the solution is known to be 
\begin{equation}\label{dufsol}
q(t) = \mathrm{sn}(\beta t, M), \qquad p(t) = \beta \mathrm{cn}(\beta t, M) \mathrm{dn}(\beta t, M),
\end{equation}
with $\mathrm{sn}$, $\mathrm{cn}$, $\mathrm{dn}$ the elliptic Jacobi functions with elliptic modulus $M=\kappa^2/\beta^2$. We consider here the parameters:
\begin{equation}\label{dufpar}
\kappa = 7, \qquad \beta = 500.
\end{equation}
At the best of our knowledge, so far such large values of the two parameters have never been considered for benchmarking. As a matter of fact, the chosen values are 100 times larger than those used, e.g., in \cite{LiWu2016} (i.e., $\kappa = 0.07$, $\beta = 5$).
For solving such problem, we compare the methods listed above by performing $N$ integration steps with a constant stepsize $h=20/N$. In Table~\ref{SVGDtab} we list the obtained result for the first four methods (i.e., SV, Gautschi, Deuflhard, and Expode), for increasing values of $N$, in terms of:
\begin{itemize}
\item execution time;
\item maximum error on $q$ ($e_q$);
\item maximum error on $p$ ($e_p$);
\item maximum error on the Hamiltonian ($e_H$).
\end{itemize}
From the listed results, one verifies that all methods have an execution time proportional to the number of steps, $N$, and very similar for the first three methods. Moreover, as expected, the accuracy of the St\"ormer-Verlet method is less than that of the Gautschi method which, in turn, is less than that of the Deuflhard method, even though the latter method evidently suffers from cancellation errors, as the stepsize is decreased. One also verifies that,  for such methods the condition $\omega h<1$ is satisfied. This latter condition is not required by Expode which, however, though 4-th order accurate, appears to be not competitive w.r.t. the Deuflhard method.

Next, in Table~\ref{Gausstab}, we list the obtained results for the $s$-stage Gauss method, $s=1,\dots,4$. They have been implemented through the SHBVM$(s,s,s)$ method. From the listed results, one verifies that, for all methods, the higher the order of the method, the smaller the execution time. Moreover, the higher-order methods are also competitive w.r.t. the previous methods. In particular, Expode appears to be not competitive w.r.t. the higher-order Gauss methods. For this reason, by also considering that its implementation is relatively difficult for higher dimensional problems, we shall not consider such method further.

\begin{figure}[t]
\centerline{\includegraphics[width=12cm,height=8cm]{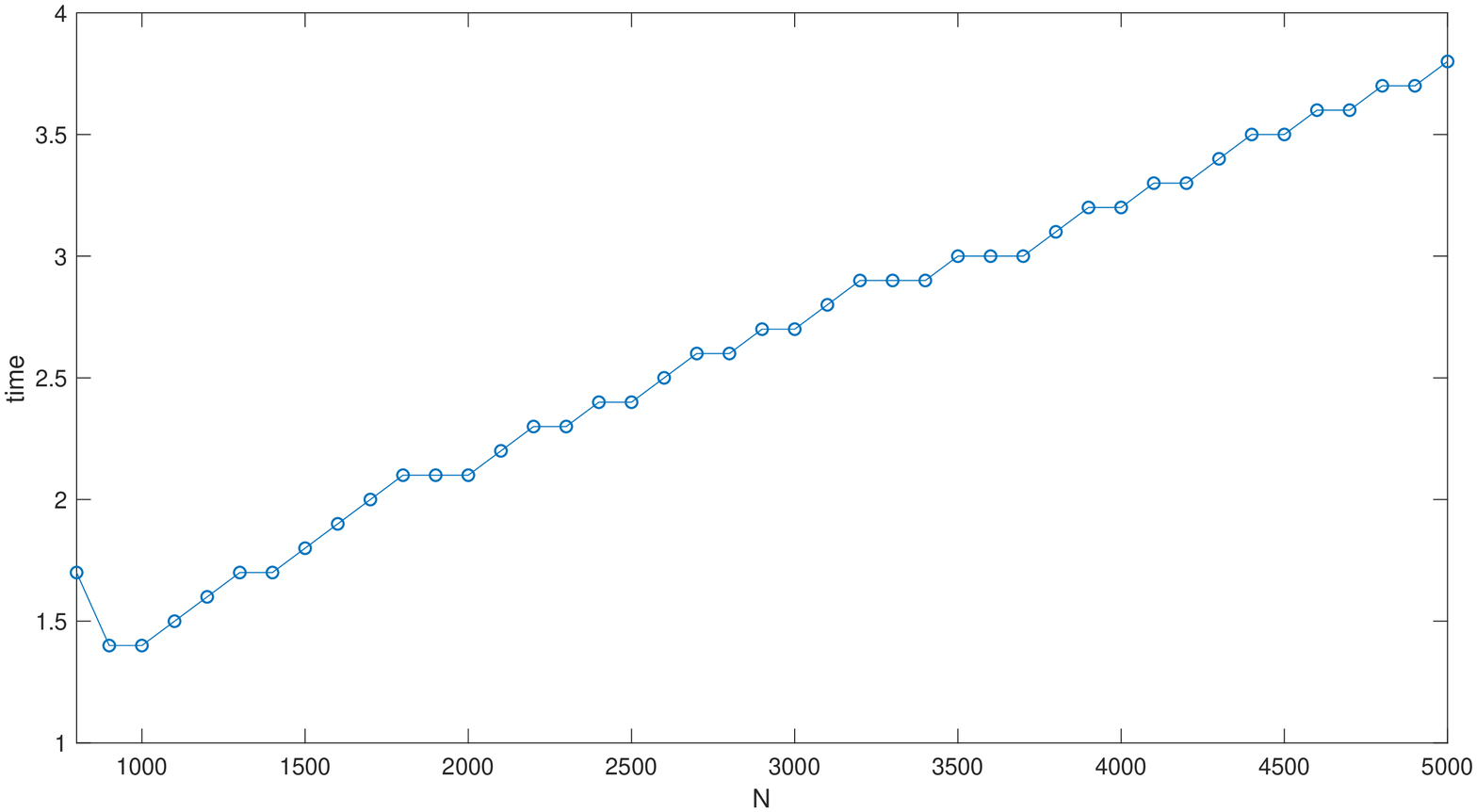}}
\caption{Problem (\ref{dufeq})  solved by the SHBVM method with stepsize $h=20/N$. Execution time versus $N$.}
\label{DufSpectralfig}
\end{figure}

At last,  in Table~\ref{Spectraltab} we list the obtained results by using the spectral method SHBVM$(k,s,s_0)$, where $(s_0, s,k)$ have been computed according to (\ref{esse}), (\ref{esse1}), and (\ref{kappa}), respectively, by considering \,$\omega = \sqrt{\kappa^2+\beta^2}$\, and \,$\nu=3$.\, 
As one may see, all errors are very small and almost constant, according to the analysis made in Sections~\ref{spectral}--\ref{HBVMs} and to what observed in Remark~\ref{norder} (in particular, the Hamiltonian error is within the round-off error level).
Moreover, also the execution times, which are very small, are almost equal, though there is a small positive trend, when reducing the stepsize $h=20/N$ (i.e., when increasing $N$). In order to make  this latter statement more precise, we computed the numerical solution by using the stepsizes
$$h = \frac{20}N, \qquad N = 800, 900, 1000, \dots, 5000.$$
In so doing, though the error remains approximately constant when increasing $N$, nevertheless, the execution time increases, as is shown in Figure~\ref{DufSpectralfig}.\footnote{Note in Figure~\ref{DufSpectralfig} that doubling the number $N$ of steps does not imply doubling the CPU time. This is due to the fact that when the number of steps is doubled, the number of stages decreases, since $\omega h$ is halved.} From that figure, one has that a value $N\approx 1000$ (to which corresponds a value $\omega h\approx 10$) seems to be the most efficient, when using SHBVM$(k,s,s_0)$. Clearly, such a method is the most efficient, among those here considered.

\begin{table}[t]
\caption{Numerical results when solving problem (\ref{dufeq}) by using a stepsize $h=20/N$.}
\label{SVGDtab}

\smallskip
\centerline{
\begin{tabular}{|r|r|r|r|r|r|r|r|r|}
\hline
\hline
\multicolumn{8}{|c|}{St\"ormer-Verlet method}\\
\hline
$N$    & time  & $e_q$  &     rate &   $e_p$   &    rate &   $e_H$  &     rate \\
\hline
1.25e6 &  18.4 & 2.65e-02 &--- & 1.30e\,01 &--- & 8.00e-05 &--- \\\
2.5e6 &  35.1 & 6.63e-03 &2.0 & 3.24e\,00 &2.0 & 2.00e-05 &2.0 \\ 
5e6 &  72.3 & 1.66e-03 &2.0 & 8.11e-01 &2.0 & 5.00e-06 &2.0 \\ 
1e7 & 143.5 & 4.15e-04 &2.0 & 2.04e-01 &2.0 & 1.25e-06 &2.0 \\ 
2e7 & 274.5 & 1.04e-04 &2.0 & 5.10e-02 &2.0 & 3.13e-07 &2.0 \\ 
\hline
\hline
\multicolumn{8}{|c|}{Gautschi method}\\
\hline
$N$    & time  & $e_q$  &     rate &   $e_p$   &    rate &   $e_H$  &     rate \\
\hline
1.25e6 &  17.2 & 7.83e-06 &--- & 3.83e-03 &--- & 7.17e-09 &--- \\ 
2.5e6 &  34.3 & 1.96e-06 &2.0 & 9.57e-04 &2.0 & 1.83e-09 &2.0 \\ 
5e6 &  68.2 & 4.89e-07 &2.0 & 2.39e-04 &2.0 & 5.43e-10 &1.8 \\ 
1e7 & 137.2 & 1.22e-07 &2.0 & 6.02e-05 &2.0 & 3.91e-10 &** \\ 
2e7 & 275.8 & 3.05e-08 &2.0 & 1.51e-05 &2.0 & 4.65e-10 &** \\
\hline
\hline
\multicolumn{8}{|c|}{Deuflhard method}\\
\hline
$N$    & time  & $e_q$  &     rate &   $e_p$   &    rate &   $e_H$  &     rate \\
\hline
6.25e5 &    8.2 & 3.26e-08 &--- & 1.60e-05 &--- & 1.67e-08 & ---\\
1.25e6 &  16.5 & 1.09e-09 &4.9 & 4.39e-07 &5.2 & 4.20e-09 &2.0 \ \\ 
2.5e6 &  32.4 & 5.12e-09 & ** & 2.50e-06 & ** & 1.05e-09 &2.0 \ \\ 
5e6 &  65.2 & 1.21e-07 & ** & 5.91e-05 & ** & 3.07e-10 &1.8 \ \\ 
1e7 & 129.8 & 5.00e-07 & ** & 2.46e-04 & ** & 8.03e-10 &** \ \\ 
\hline
\hline
\multicolumn{8}{|c|}{Expode method}\\
\hline
$N$    & time  & $e_q$  &     rate &   $e_p$   &    rate &   $e_H$  &     rate \\
\hline
6.25e3 & 5.4 & 8.58e-03 & --- & 2.81e\,00 & --- & 7.43e-04 & --- \\
1.25e4 &  10.5 & 6.07e-04 & 3.8 &1.98e-01 & 3.8 &1.68e-06 & 8.8 \\
2.5e4 &  20.8 & 3.85e-05 & 4.0 &1.25e-02 & 4.0 &8.15e-08 & 4.4 \\ 
5e4 &  45.0 & 2.43e-06 & 4.0 &8.16e-04 & 3.9 &4.24e-09 & 4.3 \\ 
1e5 &  82.8 & 1.63e-07 & 3.9 &7.67e-05 & 3.4 &9.60e-10 & ** \\ 
2e5 & 165.0 & 2.61e-08 & ** &1.22e-05 & ** &8.30e-10 & ** \\ 
\hline
\hline
\end{tabular}}
\end{table}

\begin{table}[t]
\caption{Numerical results when solving problem (\ref{dufeq}) by using a stepsize $h=20/N$.}
\label{Gausstab}

\smallskip
\centerline{
\begin{tabular}{|r|r|r|r|r|r|r|r|}
\hline
\hline
\multicolumn{8}{|c|}{1-stage Gauss method} \\
\hline
$N$  &  time &  $e_q$  &  rate &   $e_p$ &   rate &  $e_H$ &  rate\\
\hline
1.25e6 & 192.6 & 5.32e-02 & --- &2.60e\,01 & --- &3.14e-09 &---\\ 
2.5e6 & 365.6 & 1.33e-02 & 2.0 &6.51e\,00 & 2.0 &7.84e-10 & 2.0\\ 
5e6 & 700.2 & 3.33e-03 & 2.0 &1.63e\,00 & 2.0 &1.96e-10 &2.0\\ 
1e7 & 1378.9 & 8.32e-04 & 2.0 &4.09e-01 & 2.0 &4.93e-11 &2.0\\ 
2e7 & 2739.2 & 2.08e-04 & 2.0 &1.02e-01 & 2.0 &1.24e-11 &2.0\\
\hline
\hline
\multicolumn{8}{|c|}{2-stage Gauss method} \\
\hline
$N$  &  time &  $e_q$  &  rate &   $e_p$ &   rate &  $e_H$ &  rate\\
\hline
2e5 &  62.9 & 8.63e-05 & --- &4.08e-02 & --- &2.72e-11 &--- \\ 
4e5 & 114.0 & 5.40e-06 & 4.0 &2.59e-03 & 4.0 &9.35e-12 & **\\ 
8e5 & 210.1 & 3.39e-07 & 4.0 &1.65e-04 & 4.0 &2.36e-12 &**\\ 
16e5 & 403.8 & 2.12e-08 & 4.0 &1.03e-05 & 4.0 &5.44e-13 &**\\
\hline
\hline
\multicolumn{8}{|c|}{3-stage Gauss method} \\
\hline
$N$  &   time &  $e_q$  &  rate &   $e_p$ &   rate &  $e_H$ &  rate\\
\hline
2.5e4 &  11.8 & 3.98e-04 & --- &1.29e-01 & --- &4.09e-10 & ---\\
5e4 &  20.1 & 6.27e-06 & 6.0 &2.11e-03 & 5.9 &6.08e-12 & 6.1\\
1e5 &  35.1 & 9.85e-08 & 6.0 &4.66e-05 & 5.5 &1.51e-10 & **\\ 
2e5 &  65.4 & 1.47e-09 & 6.1 &6.97e-07 & 6.1 &1.09e-10 & **\\ 
\hline
\hline
\multicolumn{8}{|c|}{4-stage Gauss method} \\
\hline
$N$  &    time &  $e_q$  &  rate &   $e_p$ &   rate &  $e_H$ &  rate \\
\hline
1.25e4 &   7.3 & 6.35e-05 & --- &2.07e-02 & --- &4.68e-10 & ---\\ 
2.5e4 &  11.6 & 2.53e-07 & 8.0 &8.20e-05 & 8.0 &9.05e-13 &9.0 \\ 
5e4 &  20.3 & 9.93e-10 & 8.0 &3.34e-07 & 7.9 &2.36e-13 &** \\
1e5 &  36.1 & 1.60e-10 & ** &7.38e-08 & ** &2.16e-10 & **\\
\hline
\hline
\end{tabular}}
\end{table}

\begin{table}[t]
\caption{Numerical results when solving problem (\ref{dufeq}) by using a stepsize $h=20/N$.}
\label{Spectraltab}

\smallskip
\centerline{
\begin{tabular}{|r|r|r|r|r|r|r|}
\hline
\hline
\multicolumn{6}{|c|}{SHBVM$(k,s,s_0)$ method}\\
\hline
$N$    &  time &  $e_q$  &   $e_p$  &  $e_H$  &  $(s_0,s,k)$\\ 
\hline
   800 &  1.7 & 3.96e-10 & 7.70e-08 & 4.44e-16 & (29,50,52) \\ 
  900 &  1.4 & 5.47e-11 & 1.20e-08 & 2.22e-16 & (28,47,49) \\ 
 1000 &  1.4 & 2.70e-11 & 1.28e-09 & 4.44e-16 & (26,44,46) \\  
 1100 &  1.5 & 5.90e-11 & 2.35e-08 & 2.22e-16 & (25,42,44) \\  
 1200 &  1.6 & 1.08e-11 & 1.63e-09 & 3.33e-16 & (25,40,42) \\ 
 1300 &  1.7 & 2.63e-11 & 5.07e-09 & 4.44e-16 & (24,39,41) \\ 
 1400 &  1.7 & 2.41e-11 & 2.50e-09 & 4.44e-16 & (23,37,39) \\ 
 1500 &  1.8 & 1.77e-11 & 6.40e-09 & 4.44e-16 & (22,36,38) \\ 
\hline
\hline
\end{tabular}}
\end{table}

\subsubsection*{Fermi-Pasta-Ulam problem}
The second test problem is the well-known Fermi-Pasta-Ulam problem (see, e.g., \cite{LIMbook2016,GNI2006}), which models a physical system composed by $2m$ unit point masses disposed in series along a line, chained together by alternating weak nonlinear springs and stiff linear springs. In particular, we assume that the force exerted by the nonlinear springs is proportional to the cube of the displacement of the associated masses (cubic springs). The endpoints of the external springs are taken fixed. We denote by $q_1, q_2, .\dots , q_{2m}$ the displacements of the masses from their rest points and define the conjugate momenta as $p_i = \dot q_i$, $i = 1, . . . , 2m$. The resulting problem is Hamiltonian and is defined by the energy function
\begin{equation}\label{fpuH}
H(q,p) = \frac{1}2 \sum_{i=1}^m( p_{2i-1}^2 + p_{2i}^2) +\frac{1}2\sum_{i=1}^m\omega_i^2(q_{2i}-q_{2i-1})^2 + \sum_{i=0}^m(q_{2i+1}-q_{2i})^4, 
\end{equation}
with $q_0=q_{2m+1}=0$, and with the coefficients $\omega_i$, ruling the stiffness of the linear strings, which may be large, thus yielding a stiff oscillatory problem. We consider the parameters
\begin{equation}\label{fpupar}
m = 8, \qquad \omega_i = 10^{i-1},  \qquad \omega_{4+i} = (\pi-4+i) {\blue\cdot} 10^{4-i}, \qquad i=1,2,3,4,
\end{equation}
with the initial conditions
\begin{equation}\label{fpuy0}
q_i = \frac{i-1}{2(2m-1)}, \qquad p_i=0, \qquad i=1,\dots,2m,
\end{equation}
which evidently provide a (severe) multi-frequency highly-oscillatory problem.
A reference solution at $T=10$ has been computed by using a high-order Gauss method with a suitably small time-step.
All methods have been used to perform $N$ integration steps, with a constant stepsize $h=10/N$.

In Table~\ref{SVGDtab1} we list the obtained results for the St\"ormer-Verlet, Gautschi, and Deuflhard methods, which, as one may see, require very small stepsizes. Moreover, as for the previous problem, the Deuflhard method suffers from cancelation errors, for the smallest stepsizes used.

In Table~\ref{Gausstab1}, we list the obtained results for the $s$-stage Gauss methods, $s=1,\dots,4$. As in the previous example, only the higher order methods are relatively efficient.

\begin{figure}[t]
\centerline{\includegraphics[width=12cm,height=9cm]{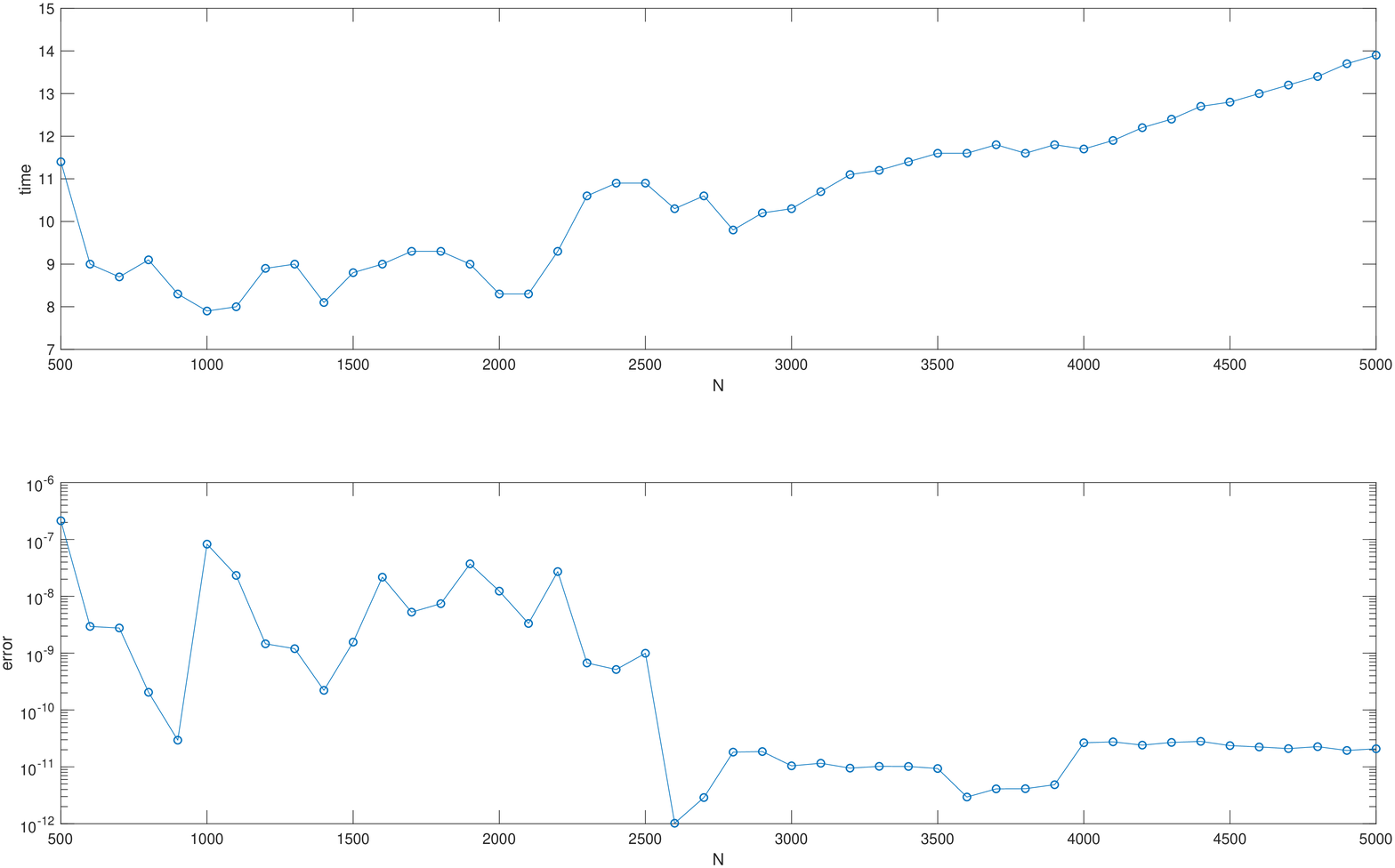}}
\caption{ Problem (\ref{fpuH})--(\ref{fpuy0})  solved by the spectral method with stepsize $h=10/N$. Upper plot: execution time versus $N$. Lower-plot: solution error versus $N$.}
\label{FPUfig}
\end{figure}

At last, in Table~\ref{Spectraltab1}, we list the obtained results by using the SHBVM$(k,s,s_0)$ method,  where $(s_0, s,k)$ have been computed according to (\ref{esse}), (\ref{esse1}), and (\ref{kappa}), respectively, by considering \,$\omega = 10^3$\, and \,$\nu=3$.\, From the listed results, one deduces that the Hamiltonian error is within the round-off error level, and the solution error is always very small, independently of the value of $N$ used, even though it seems that the optimal value of $N$ is 900 (to which corresponds a value $\omega h\approx 11$), with the minimum solution error (among those displayed) and an almost minimum execution time (approximately 8 sec). Consequently, this SHBVM method turns out to be the most efficient, among those here considered.
This fact is confirmed by the plots in Figure~\ref{FPUfig}, showing the execution time and the solution error versus $N$ (upper and lower plot, respectively), when using a stepsize  $$h = \frac{10}N, \qquad N=500,600,700,\dots,5000.$$ In fact, even though a smaller error (of the order of $10^{-11}$) is obtained for $N>2500$, a larger execution time is required.

\begin{table}[t]
\caption{Numerical results when solving problem (\ref{fpuH})--(\ref{fpuy0}) by using a stepsize $h=10/N$.}
\label{SVGDtab1}

\smallskip
\centerline{
\begin{tabular}{|r|r|r|r|r|r|r|}
\hline
\hline
\multicolumn{6}{|c|}{St\"ormer-Verlet method}\\
\hline
$10^{-4}N$    & time  & $e_y$  &     rate &   $e_H$   &    rate  \\
\hline
 16 &   7.8 & 1.71e\,00 & --- &9.25e-04 & --- \\
 32 &  15.5 & 9.63e-01 & 0.8 &2.33e-04 & 2.0 \\ 
 64 &  30.9 & 2.60e-01 & 1.9 &5.85e-05 & 2.0 \\ 
128 &  62.7 & 6.60e-02 & 2.0 &1.46e-05 & 2.0 \\ 
256 & 125.0 & 1.66e-02 & 2.0 &3.66e-06 & 2.0 \\ 
\hline
\hline
\multicolumn{6}{|c|}{Gautschi method}\\
\hline
$10^{-4}N$    & time  & $e_y$  &     rate &   $e_H$   &    rate  \\
\hline
  1 &   0.5 & 7.75e-05 & --- &1.79e-07 & --- \\
  2 &   1.0 & 1.79e-05 & 2.1 &4.08e-08 & 2.1 \\ 
  4 &   2.0 & 4.39e-06 & 2.0 &1.03e-08 & 2.0 \\ 
  8 &   4.3 & 1.09e-06 & 2.0 &2.54e-09 & 2.0 \\ 
 16 &   7.9 & 2.73e-07 & 2.0 &6.26e-10 & 2.0 \\ 
 32 &  15.8 & 6.81e-08 & 2.0 &1.48e-10 & 2.1 \\ 
 64 &  31.5 & 1.71e-08 & 2.0 &6.58e-11 & ** \\ 
128 &  62.9 & 6.32e-09 & ** &2.55e-10 & ** \\ 
\hline
\hline
\multicolumn{6}{|c|}{Deuflhard method}\\
\hline
$10^{-4}N$    & time  & $e_y$  &     rate &   $e_H$   &    rate  \\
\hline
  1 &   0.8 & 4.93e-05 & --- &4.17e-07 & --- \\
  2 &   1.5 & 1.05e-05 & 2.2 &9.73e-08 & 2.1 \\ 
  4 &   3.1 & 2.48e-06 & 2.1 &2.38e-08 & 2.0 \\ 
  8 &   6.2 & 6.11e-07 & 2.0 &5.92e-09 & 2.0 \\ 
 16 &  11.9 & 1.52e-07 & 2.0 &1.44e-09 & 2.0 \\ 
 32 &  24.0 & 2.80e-07 & ** &5.62e-10 & ** \\ 
 64 &  47.9 & 1.55e-06 & ** &8.75e-10 & ** \\ 
128 &  96.1 & 4.07e-06 & ** &2.63e-09 & ** \\ 
\hline
\hline
\end{tabular}}
\end{table}

\begin{table}[t]
\caption{Numerical results when solving problem  (\ref{fpuH})--(\ref{fpuy0}) by using a stepsize $h=10/N$.}
\label{Gausstab1}

\smallskip
\centerline{
\begin{tabular}{|r|r|r|r|r|r|}
\hline
\hline
\multicolumn{6}{|c|}{1-stage Gauss method} \\
\hline
$10^{-4}N$ &   time &  $e_y$  &  rate &  $e_H$ & rate\\
\hline
 32 & 107.2 & 2.16e\,00 & --- &4.90e-11& --- \\
 64 & 213.5 & 5.44e-01 & 2.0 &1.27e-11 &1.9 \\ 
128 & 428.2 & 1.34e-01 & 2.0 &3.33e-12 &1.9\\ 
256 & 855.2 & 3.32e-02 & 2.0 &2.04e-12 &**\\ 
\hline
\hline
\multicolumn{6}{|c|}{2-stage Gauss method} \\
\hline
$10^{-4}N$  &   time &  $e_y$  &  rate &  $e_H$ & rate\\
\hline
  4 &  27.7 & 1.42e-01 & --- &5.43e-12 &---\\
  8 &  50.0 & 8.83e-03 & 4.0 &4.51e-13 & 3.6\\
 16 &  95.4 & 5.52e-04 & 4.0 &3.04e-13 & **\\ 
 32 & 187.6 & 3.45e-05 & 4.0 &3.89e-13 & **\\ 
 64 & 354.9 & 2.16e-06 & 4.0 &6.08e-13 & **\\ 
128 & 661.5 & 1.35e-07 & 4.0 &2.74e-13 & **\\ 
\hline
\hline
\multicolumn{6}{|c|}{3-stage Gauss method} \\
\hline
$10^{-4}N$  &  time &  $e_y$  &  rate &  $e_H$ & rate\\
\hline
  1 &  11.1 & 2.52e-01 & --- &1.35e-11 & ---\\
  2 &  18.9 & 4.00e-03 & 6.0 &2.46e-13 & 5.8\\ 
  4 &  33.0 & 6.30e-05 & 6.0 &1.33e-13 & **\\ 
  8 &  56.6 & 9.86e-07 & 6.0 &1.21e-13 & **\\ 
 16 & 106.4 & 1.54e-08 & 6.0 &1.27e-13 & **\\ 
 32 & 211.9 & 2.76e-10 & 5.8 &4.17e-13 & **\\ 
\hline
\hline
\multicolumn{6}{|c|}{4-stage Gauss method} \\
\hline
$10^{-4}N$ &   time &  $e_y$  &  rate &  $e_H$ & rate\\
 \hline
  1 &  11.8 & 9.96e-04 & --- &4.71e-13 & ---\\
  2 &  19.9 & 3.98e-06 & 8.0 &1.27e-13 & **\\ 
  4 &  34.8 & 1.56e-08 & 8.0 &8.79e-14 & **\\
  8 &  59.6 & 5.28e-11 & 8.2 &1.36e-13 & **\\ 
 16 & 109.8 & 8.00e-12 & ** &1.12e-13 & **\\ 
\hline
\hline
\end{tabular}}
\end{table}

\begin{table}[t]
\caption{Numerical results when solving problem (\ref{fpuH})--(\ref{fpuy0}) by using a stepsize $h=10/N$.}
\label{Spectraltab1}

\smallskip
\centerline{
\begin{tabular}{|r|r|r|r|r|r|r|}
\hline
\hline
\multicolumn{5}{|c|}{SHBVM$(k,s,s_0)$ method}\\
\hline
$N$    &  time &  $e_y$  &  $e_H$  &  $(s_0,s,k)$ \\  
\hline
  500 & 11.4 & 2.13e-07 & 1.78e-15 & (36,66,68) \\  
  600 &  9.0 & 2.95e-09 & 1.78e-15 & (33,59,61) \\  
  700 &  8.7 & 2.77e-09 & 1.78e-15 & (31,54,56) \\  
  800 &  9.1 & 2.05e-10 & 2.00e-15 & (29,50,52) \\   
  900 &  8.3 & 2.95e-11 & 1.78e-15 & (28,47,49) \\  
 1000 &  7.9 & 8.28e-08 & 1.78e-15 & (26,44,46) \\   
 1100 &  8.0 & 2.33e-08 & 1.78e-15 & (25,42,44) \\ 
 1200 &  8.9 & 1.46e-09 & 2.00e-15 & (25,40,42) \\ 
 1300 &  9.0 & 1.20e-09 & 1.78e-15 & (24,39,41) \\  
 1400 &  8.1 & 2.22e-10 & 1.78e-15 & (23,37,39) \\ 
 1500 &  8.8 & 1.56e-09 & 2.00e-15 & (22,36,38) \\ 
\hline
\hline
\end{tabular}}
\end{table}

\subsubsection*{Nonlinear Schr\"odinger equation}

At last, we consider a highly-oscillatory problem, in the form (\ref{first}), deriving from the space semi-discretization of the nonlinear Sch\"odinger equation \cite{BBFCI2018}, which we sketch below in a much simplified form.\footnote{This is due to the particular initial condition considered.} The problem we want to solve is
\begin{eqnarray}\label{nlse}
\ii\psi_t(x,t) +\psi_{xx}(x,t) + \kappa |\psi(x,t)|^2\psi(x,t) &=& 0, \qquad \qquad (x,t)\in[0,2\pi]\times[0,5], \\
\psi(x,0) &=& \ee^{\ii r  x},  \nonumber
\end{eqnarray}
coupled with periodic boundary conditions. The solution of this ``toy'' problem is readily seen to be given by
$$
\psi(x,t)  = \ee^{i(r  x-\mu t)}, \qquad \mu = r ^2-\kappa.
$$
In particular, in order to simplify the arguments, we consider the values:
$$
r  = 20, \qquad \kappa = \frac{\pi}{10} \qquad \Rightarrow\qquad \mu \approx 399.7.
$$
To begin with, we separate the real and imaginary parts of the solution,
$$ \psi(x,t) = u(x,t) + \ii v(x,t),$$
and consider their expansions (in space) along an orthonormal Fourier basis, 
$$
u(x,t) = \sum_{j=0}^r  \left[c_j(x) \xi_j(t)+  s_j(x)\eta_j(t)\right],\quad 
v(x,t) = \sum_{j=0}^r  \left[c_j(x) \alpha_j(t)+  s_j(x)\beta_j(t)\right],\quad x\in[0,2\pi],
$$
with
$$c_j(x) = \sqrt{\frac{2-\delta_{j0}}{2\pi}}\cos(j x), \qquad s_j(x) = \sqrt{\frac{2-\delta_{j0}}{2\pi}}\sin(j x), \qquad j=0,\dots,r ,$$
$\delta_{j0}$ the Kronecker delta, and unknown time-dependent coefficients $\xi_j(t), \eta_j(t), \alpha_j(t), \beta_j(t)$.\footnote{For this simple problem, actually only the coefficients $\xi_r (t),\eta_r (t),\alpha_r (t),\beta_r (t)$ are nonzero.}  Subsequently, by setting the vectors
$$\bfq(t) = \pmatrix{c} \xi_0(t)\\ \vdots \\ \xi_r (t)\\ \eta_1(t)\\ \vdots \\ \eta_r (t)\endpmatrix, ~
\bfp(t) = \pmatrix{c} \alpha_0(t)\\ \vdots \\ \alpha_r (t)\\ \beta_1(t)\\ \vdots \\ \beta_r (t)\endpmatrix, ~
\bfw(t) = \pmatrix{c} c_0(t)\\ \vdots \\ c_r (t)\\ s_1(t)\\ \vdots \\ s_r (t)\endpmatrix ~\in\RR^{2r +1},\qquad \bfy = \pmatrix{c} \bfq\\ \bfp\endpmatrix,$$
and the matrix
$$D = \pmatrix{cccccc} 0\\ &1\\ &&1 \\ &&&\ddots \\ &&&& r  \\ &&&&&r \endpmatrix\in\RR^{2r +1\times 2r +1},$$
we obtain, by setting $e_i\in\RR^2$ the $i$-th unit vector, $i=1,2$,
$$|\psi(x,t)|^2 = \left[\left((e_1 \otimes \bfw(x))^\top \bfy(t)\right)^2 + \left((e_2 \otimes \bfw(x))^\top \bfy(t)\right)^2 \right],$$
thus arriving at the Hamiltonian system of ODEs
\begin{equation}\label{nlse1}
\dot\bfy = J_2\otimes D^2\,\bfy - \kappa \int_0^{2\pi} J_2\otimes (\bfw(x)\bfw(x)^\top) \, \bfy  \left[\left((e_1 \otimes \bfw(x))^\top \bfy(t)\right)^2 + \left((e_2 \otimes \bfw(x))^\top \bfy(t)\right)^2 \right]\dd x,
\end{equation}
(here $J_2$ is the same matrix as that defined in (\ref{sec1})), with Hamiltonian
\begin{equation}\label{Hnlse}
H(\bfy) = \frac{1}2\left( \bfy^\top I_2\otimes D^2\,\bfy - \frac{\kappa}2\int_0^{2\pi}  \left[\left((e_1 \otimes \bfw(x))^\top \bfy(t)\right)^2 + \left((e_2 \otimes \bfw(x))^\top \bfy(t)\right)^2 \right]^2\dd x\right).
\end{equation}
This latter function, in turn, is equivalent to the Hamiltonian functional defining (\ref{nlse}), i.e.,\footnote{The reader is referred to \cite{BBFCI2018} for full details.}
$$\H[\psi] = \frac{1}2\int_0^{2\pi} |\psi_x|^2 -\frac{\kappa}2|\psi|^4\dd x.$$
At last, in order to derive a fully discrete problem, the integrals in (\ref{nlse1})--(\ref{Hnlse}) are (exactly) computed via a composite trapezoidal rule, at the abscissae
\begin{equation}\label{xim} x_\ell = \ell\frac{2\pi}m, \qquad \ell=0,\dots,m, \qquad m = 4r +1.\end{equation}

Because of the structure of the problem, to implement the SHBVM method we shall consider the parameters
$$\omega = 400 \equiv r ^2 = \|D^2\|, \qquad \nu = 1,$$
so that, in this case, the two parameters $s_0$ and $s$ (see (\ref{esse}) and (\ref{esse1}), respectively) coincide. The parameter $k$, in turn, is computed, as usual, from (\ref{kappa}).

In Table~\ref{Gausstab2}, we list the obtained results for the $s$-stage Gauss methods, $s=1,\dots,4$, by using a time-step $h=5/N$, for increasing values of $N$. As in the previous examples, only the higher order methods are relatively efficient, even though the Hamiltonian error is always very small.

At last, in Table~\ref{Spectraltab2}, we list the obtained results by using the SHBVM$(k,s,s)$ method,\footnote{In fact, in this case $s_0=s$.} by using the time-steps:
$$h = 5/N, \qquad N=200,250,300,350,400,450,500.$$ From the listed results, one deduces that the Hamiltonian error is always within the round-off error level. Moreover, also the solution error is always very small and uniform, independently of the considered  value of $N$, even though it seems that the optimal value of $N$ is 250 (to which corresponds a value $\omega h\approx 8$), with an almost minimum solution error (among those displayed) and a minimum execution time (2.7 sec). Consequently, with such a value of $N$, SHBVM(26,24,24) is clearly the most efficient method, among those here considered.

\begin{table}[p]
\caption{Numerical results when solving problem  (\ref{nlse1})--(\ref{Hnlse}) by using a stepsize $h=5/N$.}
\label{Gausstab2}

\smallskip
\centerline{
\begin{tabular}{|r|r|r|r|r|}
\hline
\hline
\multicolumn{5}{|c|}{1-stage Gauss method} \\
\hline
$10^{-3}N$ &   time &  $e_y$  &  rate &  $e_H$ \\
\hline
16 &  6.3 & 2.01e\,02 &  --- & 8.55e-15\\ 
32 & 10.9 & 6.65e\,01 &  1.6 & 1.22e-14\\ 
64 & 21.7 & 1.69e\,01 &  2.0 & 2.10e-14\\ 
128 & 44.6 & 4.23e\,00 &  2.0 & 3.85e-14\\ 
256 & 88.7 & 1.06e\,00 &  2.0 & 5.02e-14\\ 
512 & 174.4 & 2.64e-01 &  2.0 & 2.35e-14\\ 
1024 & 290.4 & 6.61e-02 &  2.0 & 1.38e-13\\ 
\hline
\hline
\multicolumn{5}{|c|}{2-stage Gauss method} \\
\hline
$10^{-3}N$  &   time &  $e_y$  &  rate &  $e_H$ \\
\hline
 4 &  4.0 & 1.77e\,01 &  --- & 1.56e-13\\ 
8 &  6.8 & 1.12e\,00 &  4.0 & 7.75e-14\\ 
16 & 11.1 & 7.02e-02 &  4.0 & 6.95e-14\\ 
32 & 19.5 & 4.39e-03 &  4.0 & 8.44e-15\\ 
64 & 34.7 & 2.75e-04 &  4.0 & 1.78e-14\\ 
128 & 68.9 & 1.72e-05 &  4.0 & 2.82e-14\\ 
256 & 120.9 & 1.07e-06 &  4.0 & 2.49e-14\\ 
\hline
\hline
\multicolumn{5}{|c|}{3-stage Gauss method} \\
\hline
$10^{-4}N$  &  time &  $e_y$  &  rate &  $e_H$ \\
\hline
2 &  3.1 & 1.97e\,00 &  --- & 2.00e-14\\ 
4 &  5.1 & 3.17e-02 &  6.0 & 2.26e-14\\ 
8 &  8.1 & 4.99e-04 &  6.0 & 1.93e-14\\ 
16 & 14.7 & 7.82e-06 &  6.0 & 7.33e-15\\ 
32 & 26.4 & 1.22e-07 &  6.0 & 7.11e-15\\ 
64 & 48.6 & 2.20e-09 &  5.8 & 1.42e-14\\ 
128 & 84.4 & 3.28e-10 &  ** & 2.24e-14\\ 
\hline
\hline
\multicolumn{5}{|c|}{4-stage Gauss method} \\
\hline
$10^{-3}N$ &   time &  $e_y$  &  rate &  $e_H$ \\
 \hline
1 &  2.1 & 1.84e\,00 &  --- & 3.97e-13\\ 
2 &  3.0 & 7.88e-03 &  7.9 & 5.48e-14\\ 
4 &  4.7 & 3.14e-05 &  8.0 & 3.87e-14\\ 
8 &  8.2 & 1.23e-07 &  8.0 & 9.99e-15\\ 
16 & 13.7 & 4.33e-10 &  8.2 & 5.22e-15\\ 
32 & 25.4 & 1.01e-10 &  ** & 1.40e-14\\ 
\hline
\hline
\end{tabular}}
\end{table}

\begin{table}[t]
\caption{Numerical results when solving problem (\ref{nlse1})--(\ref{Hnlse}) by using a stepsize $h=5/N$.}
\label{Spectraltab2}

\smallskip
\centerline{
\begin{tabular}{|r|r|r|r|r|r|r|}
\hline
\hline
\multicolumn{5}{|c|}{SHBVM$(k,s,s)$ method}\\
\hline
$N$    &  time &  $e_y$  &  $e_H$  &  $(s,k)$ \\  
\hline
 200 &  2.7 & 1.50e-10 & 4.44e-16 & (26,28)\\ 
 250 &  2.7 & 4.94e-11 & 4.44e-16 & (24,26)\\ 
 300 &  2.9 & 2.43e-10 & 4.44e-16 & (22,24)\\ 
 350 &  3.0 & 1.43e-10 & 4.44e-16 & (21,23)\\ 
 400 &  3.1 & 4.83e-11 & 3.33e-16 & (20,22)\\ 
 450 &  3.5 & 4.33e-11 & 4.44e-16 & (19,21)\\
 500 &  3.6 & 5.53e-11 & 4.44e-16 & (19,21)\\
\hline
\hline
\end{tabular}}
\end{table}

\section{Concuding remarks}\label{fine}

In this paper, we have seen that spectral methods in time can be very efficient for solving, possibly multi-frequency, highly-oscillatory Hamiltonian problems. In particular, their implementation via a slight modification of Hamiltonian Boundary Value Methods, provides a very efficient energy-conserving procedure, able to solve such problems up to the maximum accuracy allowed by the used finite-precision arithmetic. Some numerical tests on some severe tests problems duly confirm the theoretical achievements. 

It is worth noticing that the proposed approach could be in principle used also for solving general ODE-IVPs, thus providing a spectrally accurate method of solution in time. This will be the subject of future investigations.

\medskip
{\bf Acknowledgements.} This paper emerged during visits of the first author at the Departamento  Matem\'{a}tica Aplicada, Universidad de Zaragoza, in July and October 2017. The authors wish to thank Gianmarco Gurioli, for carefully reading the manuscript.


\begin{thebibliography}{99}
\setlength{\itemsep}{0.05cm}

\bibitem{ABI2015} P.\,Amodio, L.\,Brugnano, F.\,Iavernaro. Energy-conserving methods for Hamiltonian Boundary Value Problems and applications in astrodynamics. {\em  Adv. Comput. Math.} {\bf 41} (2015) 881--905.

\bibitem{ArEnKiLeTs13} G.\,Ariel, B.\,Engquist, S.\,Kim, Y.\,Lee, R.\,Tsai.  A Multiscale Method for Highly Oscillatory Dynamical Systems Using a Poincar\'{e} Map Type Technique, {\em J. Sci. Comput.}, {\bf 54} (2013) 247--268.

\bibitem{BBFCI2018} L.\,Barletti, L.\,Brugnano, G.\,Frasca Caccia, F.\,Iavernaro. Energy-conserving methods for the nonlinear Schr\"odinger equation. {\em Appl. Math.Comput.}  {\bf 318} (2018) 3--18.
 
\bibitem{BeSt2000} P.\,Betsch, P.\,Steinmann. Inherently Energy Conserving Time Finite Elements for Classical Mechanics. {\em J. Comp. Phys.} {\bf 160} (2000) 88--116.

\bibitem{Bo1997} C.L.\,Bottasso. A new look at finite elements in time: a variational interpretation of Runge--Kutta methods. {\em Appl. Numer. Math.} {\bf 25} (1997) 355--368.

\bibitem{Br2000} L.\,Brugnano. Blended Block BVMs (B$_3$VMs): A Family of Economical Implicit Methods for ODEs. {\em J. Comput. Appl. Math.} {\bf  116} (2000) 41--62.

\bibitem{BCMR2012} L.\,Brugnano, M.\,Calvo, J.I.\,Montijano, L.\,R\`andez.  Energy preserving methods for Poisson systems. {\em J. Comput. Appl. Math.} {\bf  236} (2012) 3890--3904.

\bibitem{BCMR16} L.\,Brugnano, M.\, Calvo, J. I.\, Montijano, L.\,R\'{a}ndez.  Fourier Methods for Oscillatory Differential Problems with a Constant High Frequency. {\em AIP Conference Proc.} {\bf 1863} (2017) 020003.

\bibitem{BFCI2014} L.\,Brugnano, G.\,Frasca Caccia, F.\,Iavernaro. Efficient  implementation of Gauss collocation and Hamiltonian Boundary Value Methods. {\em Numer. Algorithms} {\bf 65} (2014) 633--650. 

\bibitem{BFCI2015} L.\,Brugnano, G.\,Frasca Caccia, F.\,Iavernaro.  Energy conservation issues in the numerical solution of the semilinear wave equation. {\em Appl. Math.Comput.} {\bf 270} (2015) 842--870

\bibitem{BrIa2012} L.\,Brugnano, F.\,Iavernaro.  Line Integral Methods which preserve all invariants of conservative problems.  {\em J. Comput. Appl. Math.} {\bf  236} (2012) 3905--3919.

\bibitem{LIMbook2016} L.\,Brugnano, F.\,Iavernaro. {\em Line Integral Methods for Conservative Problems}.  Chapman et Hall/CRC, Boca Raton, FL, 2016. 

\bibitem{BIT2010} L.\,Brugnano, F.\,Iavernaro, D.\,Trigiante.  Hamiltonian Boundary Value Methods (Energy Preserving Discrete Line Integral Methods). {\em  JNAIAM. J. Numer. Anal. Ind. Appl. Math.} {\bf 5}, No.\,1-2 (2010)  17--37.

\bibitem{BIT2011} L.\,Brugnano, F.\,Iavernaro, D.\,Trigiante. A note on the efficient implementation of Hamiltonian BVMs. {\em J. Comput. Appl. Math.} {\bf 236} (2011) 375--383.

\bibitem{BIT2012_0} L.\,Brugnano, F.\,Iavernaro, D.\,Trigiante. The lack of continuity and the role of Infinite and infinitesimal in numerical methods for ODEs: the case of symplecticity.  {\em Appl. Math.Comput.} {\bf 218} (2012)  8053--8063.

\bibitem{BIT2012} L.\,Brugnano, F.\,Iavernaro, D.\,Trigiante.  A simple framework for the derivation and analysis of effective one-step methods for ODEs. {\em Appl. Math.Comput.} {\bf 218} (2012) 8475--8485.

\bibitem{BIT2012_1} L.\,Brugnano, F.\,Iavernaro, D.\,Trigiante. A two-step, fourth-order method with energy preserving properties. {\em  Comput. Phys. Commun.} {\bf 183} (2012) 1860--1868.

\bibitem{BIT2012_2} L.\,Brugnano, F.\,Iavernaro, D.\,Trigiante. Energy and QUadratic Invariants Preserving integrators based upon Gauss collocation formulae. {\em SIAM J. Numer. Anal.} {\bf 50}, No.\,6 (2012) 2897--2916.

\bibitem{BIT2015} L.\,Brugnano, F.\,Iavernaro, D.\,Trigiante. Analisys of Hamiltonian Boundary Value Methods (HBVMs): a class of energy-preserving Runge-Kutta methods for the numerical solution of polynomial Hamiltonian systems.  {\em Commun. Nonlinear Sci. Numer. Simul.} {\bf 20} (2015) 650--667.

\bibitem{BrMa2002} L.\,Brugnano, C.\,Magherini. Blended Implementation of Block Implicit Methods for ODEs. {\em Appl. Numer. Math.} {\bf 42} (2002) 29--45.

\bibitem{BrMa2004} L.\,Brugnano, C.\,Magherini. The BiM code for the numerical solution of ODEs. {\em J. Comput. Appl. Math.} {\bf  164-165} (2004) 145--158.

\bibitem{BrMa2007} L.\,Brugnano, C.\,Magherini. Blended Implicit Methods for solving ODE and DAE problems, and their extension for second order problems. {\em J. Comput. Appl. Math.} {\bf  205} (2007) 777--790.

\bibitem{BrMa2009} L.\,Brugnano, C.\,Magherini. Recent advances in linear analysis of convergence for splittings for solving ODE problems. {\em Appl. Numer. Math.} {\bf 59} (2009) 542--557.

\bibitem{BrMaMu2006} L.\,Brugnano, C.\,Magherini, F.\,Mugnai. Blended Implicit Methods for the Numerical Solution of DAE Problems. {\em J. Comput. Appl. Math.} {\bf 189} (2006) 34--50.

\bibitem{CaChMuSS11} M.P.\,Calvo, P.\,Chartier, A.\,Murua, J.M.\,Sanz-Serna. Numerical stroboscopic averaging for ODEs and DAEs. {\em Appl. Numer. Math. }, {\bf 61}, 10 (2011) 1077--1095.

\bibitem{CaJaMoRa04} M.\,Calvo, L.O.\,Jay, J.I.\,Montijano, L.\,R\'{a}ndez. Approximate compositions of a near identity map by multi-revolution Runge-Kutta methods. {\em Numer. Math.} {\bf 4} (2004) 635--666.

\bibitem{CaGo15} B.\,Cano, A.\,Gonz\'{a}lez-Pach\'{o}n. Exponential time integration of solitary waves of cubic Schr\"{o}dinger equation. {\em Appl. Numer. Math.}, {\bf 91} (2015) 26--45.

\bibitem{Cano2013}  B.\,Cano, M.J.\,Moreta. 
High-order symmetric multistep cosine methods.  {\em Appl. Numer. Math.} {\bf 66} (2013) 30--44. 

\bibitem{CaChFa09} F.\,Castella, P.\,Chartier, E.\,Fao. An averaging technique for highly oscillatory Hamiltonian problems. {\em SIAM J. Numer. Math.} {\bf 47}, 4 (2009) 2808--2837.

\bibitem{Cohen2006} D.\,Cohen, T.\,Jahnke, K.\,Lorenz, C.\,Lubich. Numerical Integrators for Highly Oscillatory Hamiltonian Systems: A Review. In: {\em Analysis, Modeling and Simulation of Multiscale Problems}, A.\,Mielke (eds). Springer, Berlin, Heidelberg (2006) pp.\,553--576.

\bibitem{De1979} P.\,Deuflhard. A study of extrapolation methods based on multistep schemes without parasitic solutions. {\em Z. angew. Math. Phys.} {\bf 30} (1979) 177--189.

\bibitem{EW1999} G.A.\,Evans, J.R.\,Webster. A comparison of some methods for the evaluation of highly oscillatory integrals. {\em J. Comput. Appl. Math.} {\bf 112} (1999) 55--69.

\bibitem{Ga1961} W.\,Gautschi. Numerical integration of ordinary differential equations based on trigonometric polynomials. {\em Numer. Math.} {\bf 3} (1961) 381--397.

\bibitem{GaSaSk98} B.\,Garc\'{\i}a-Archilla, J.M.\,Sanz-Serna, R.D.\,Skeel.  Long-time steps methods for oscillatory differential equations.
    {\em SIAM J. Sci. Comput.},{\bf  20}, 3 (1998) 930--963.

\bibitem{Ha2010} E.\,Hairer. Energy-Preserving Variant of Collocation Methods.   {\em  JNAIAM. J. Numer. Anal. Ind. Appl. Math.} {\bf 5}, No.\,1-2 (2010)  73--84.

\bibitem{GNI2006} E.\,Hairer, C.\,Lubich, G.\,Wanner. {\em Geometric Numerical Integration. Structure-Preserving Algorithms for Ordinary Differential Equations, Second ed.}, Springer, Berlin, 2006.

\bibitem{HoOs2005} M.\,Hochbruck, A.\,Ostermann. Explicit exponential Runge-Kutta methods for semilinear
parabolic problems. {\em SIAM J. Numer. Anal.} {\bf 43} (2005) 1069--1090.

\bibitem{HoOs10} M.\,Hochbruck, A.\,Ostermann.   Exponential integrators.    {\em Acta Numer.} {\bf  19} (2010) 209--286.

\bibitem{Hu1972} B.L.\,Hulme. One-Step Piecewise Polynomial Galerkin Methods for Initial Value Problems. {\em Math. Comp.} {\bf 26} (1972) 415--426.

\bibitem{Hu1972-1} B.L.\,Hulme. Discrete Galerkin and related one-step methods for ordinary differential equations. {\em Math. Comp.} {\bf 26} (1972) 881--891.

\bibitem{Jan2014} G.\,Jansing. EXPODE -- Advanced Exponential Time Integration Toolbox for MATLAB. {\tt arXiv:1404.4580v1 [math.NA]}, 2014

\bibitem{LiWu2016} Yu-wen Li, Xinyuan Wu. Functionally fitted energy--preserving methods for solving oscillatory nonlinear Hamiltonian systems. {\em SIAM J.\,Numer.\,Anal.} {\bf 54}, No.\,4 (2016) 2036--2059.

\bibitem{PeJaYe97} L.R.\,Petzold, L.O.\,Jay, J.\,Yen. Numerical solution of highly oscillatory ordinary differential equations. {\em Acta Numer.} {\bf 6} (1997) 437--483.

\bibitem{SSC1994} J.M.\,Sanz-Serna, M.P.\,Calvo. {\em Numerical Hamiltonian problems.} Chapman \& Hall, London, 1994.

\bibitem{ShHaJa04} S.S.\,Shome, E.J.\,Haug, L.O.\,Jay. Dual-rate integration using partitioned runge-kutta methods for mechanical systems with interacting subsystems. {\em Mechanics Based Design of Structures and Machines}, {\bf 32}, 3 (2004) 253--282.

\bibitem{SoMa2011} T.\,Solcia, P.\,Masarati. Multirate simulation of complex multibody systems. In ECCOMAS Thematic Conference, 2011.

\bibitem{TaSu2012} W.\,Tang, Y.\,Sun. Time finite element methods: a unified framework for numerical discretizations of ODEs. {\em Appl. Math. Comp.} {\bf 219} (2012) 2158--2179.

\bibitem{TuBe91} M.E.\,Tuckerman, B.J.\,Berne. Molecular dynamics in systems with multiple time scales: Systems with stiff and soft degrees of freedom and with short and long range forces. {\em J. Chem. Phys.} {\bf 95}, 11 (1991) 8362--8364.
\end{thebibliography}
\end{document}